\documentclass[11pt]{article}
\usepackage[a4paper,margin=2.5cm]{geometry}
\usepackage[T1]{fontenc}
\usepackage{lmodern,microtype,amsthm}
\usepackage{amsmath,amssymb,mathtools}
\usepackage{graphicx,tabularx,booktabs,multirow,pgfplotstable,tikz,tikz-cd,enumitem}
\usetikzlibrary{arrows.meta,positioning,calc,fit}
\pgfplotsset{compat=1.18}
\setlist{label=(\alph*)}
\def\datapath{}\def\figwidth{\linewidth}
\usepackage{natbib,placeins}


\bibpunct[, ]{(}{)}{,}{a}{}{,}

\usepackage[colorlinks=true,allcolors=blue,hypertexnames=false]{hyperref}
\definecolor{mycolor1}{rgb}{0.105882,0.619608,0.466667}
\definecolor{mycolor2}{rgb}{0.85098,0.372549,0.00784314}
\definecolor{mycolor3}{rgb}{0.458824,0.439216,0.701961}
\definecolor{mycolor4}{rgb}{0.905882,0.160784,0.541176}
\definecolor{mycolor5}{rgb}{0.4,0.65098,0.117647}
\definecolor{mycolor6}{rgb}{0.65098,0.462745,0.113725}
\definecolor{mycolor7}{rgb}{0.901961,0.670588,0.00784314}
\definecolor{mycolor8}{rgb}{0.4,0.4,0.4}
\definecolor{mycolor9}{rgb}{0.301961,0,0.294118}
\definecolor{mycolor10}{rgb}{0.0313725,0.25098,0.505882}
\newtheorem{theorem}{Theorem}
\newtheorem{lemma}{Lemma}
\newtheorem{proposition}{Proposition}
\newtheorem{corollary}{Corollary}

\theoremstyle{definition}

\newtheorem{remark}{Remark}

\newenvironment{namedproof}[1]{\begin{proof}[#1]}{\end{proof}}
\title{Preconditioning for Quadratic Matrix Optimization via Constraint Reshaping}
\author{Zhen Peng\thanks{State Key Laboratory of Mathematical Sciences, Academy of Mathematics and Systems Science, Chinese Academy of Sciences, China (\href{mailto:pengzhen@lsec.cc.ac.cn}{pengzhen@lsec.cc.ac.cn}, \href{mailto:gaobin@lsec.cc.ac.cn}{gaobin@lsec.cc.ac.cn}).}\and Bin Gao\footnotemark[1]}
\hypersetup{pdftitle={Preconditioning for Quadratic Matrix Optimization via Constraint Reshaping},pdfauthor={Zhen Peng and Bin Gao}}
\date{}

\newcommand{\prooflocation}{appendices}
\begin{document}
\maketitle
\begin{abstract}
We study quadratic optimization over orthogonality or symplectic constraints, a problem arising in eigenvalue computation, subspace analysis, and reduced-order modeling. 
    Ill-conditioned problems have a wide positive Hessian spectrum at a minimizer, which slows first-order Riemannian methods. 
    Existing preconditioned-metric approaches choose a Riemannian metric, which determines the resulting Hessian spectrum. 
    We instead redesign the problem itself to find the narrowest possible spectrum. 
    We reshape the constraints into equivalent formulations that preserve the optimal value and the solution set, giving a family of preconditioned problems parametrized by positive constraint weights. 
    We derive a closed form for the positive Hessian spectrum in terms of these weights. 
    Minimizing the condition number then reduces to a linear program (LP) in the logarithms of the weights. 
    The LP attains a globally optimal design within the paired equal-block family, and its strong duality gives matching lower bounds on the condition number. 
    On tested instances, the optimized formulation substantially reduces the local condition number, with returned spectra matching predictions on both orthogonal and symplectic instances. 
    We also benchmark the design against Lagrangian and left-right preconditioned metrics. 
    The method thus provides an exact, optimization-based alternative to heuristic preconditioned metrics.
\end{abstract}
\section{Introduction}
\label{sec:introduction}
    We study the quadratic trace problem with orthogonality or symplectic constraints
    \begin{equation}\label{eq:quadratic-problem}\tag{St/Sp}
        \min_Y h(Y)=\mathrm{tr}(Y^\top CYD),\quad \mathrm{s.}\ \mathrm{t.}\quad Y\in\mathrm{St}(p,n)\ \text{or}\ Y\in\mathrm{Sp}(p,n),
    \end{equation}
    where $C\in\mathbb R^{n\times n}$ and $D\in\mathbb R^{p\times p}$ are symmetric positive definite (SPD), with $1\le p\le n$.
    The manifolds are $\mathrm{St}(p,n)=\{Y:Y^\top Y=I_p\}$ and $\mathrm{Sp}(p,n)=\{Y:Y^\top J_nY=J_p\}$, with $J_m=I_{m/2}\otimes\left[\begin{smallmatrix}0&1\\-1&0\end{smallmatrix}\right]$ and even $n,p$ in the symplectic case.
    The model recovers symmetric and symplectic eigenvalue problems when $D=I$, Brockett costs for diagonal $D$, and more general quadratic trace problems~\citep{edelman1998geometry,son2021computing,liu2019quadratic,son2025brockett}.
    Beyond eigenvalue computation, these problems appear in subspace analysis and reduced-order modeling~\citep{edelman1998geometry,gao2021riemannian}.
    At a minimizer, the positive Hessian eigenvalues depend on the ratios and pairwise sums and differences of the eigenvalues of $C$ and $D$. 
    A wide eigenvalue spread enlarges the largest positive Hessian eigenvalue, while small pairwise gaps shrink the smallest one. 
    Both effects raise the positive Hessian condition number, which is the ill-conditioning addressed in this paper.

\begin{figure}[htbp]
        \centering
        \includegraphics[width=0.95\linewidth]{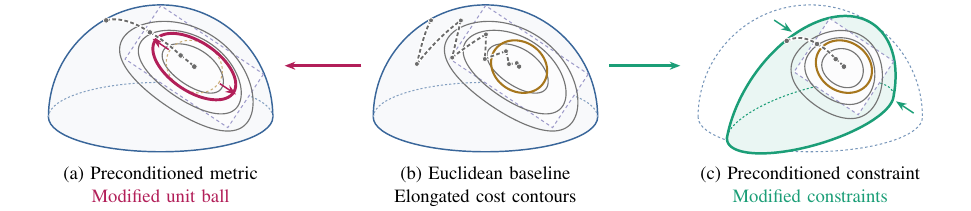}
        \caption{Schematic local effects of preconditioning. (b) The original ill-conditioned problem, with elongated cost contours and a slowly converging zigzag trajectory. (a) The feasible set is fixed and the Riemannian metric is changed, reshaping the local unit ball and straightening the trajectory. (c) The constraints and the corresponding cost contours are reshaped, giving a more direct trajectory in the equivalent formulation.}
        \label{fig:landscape-routes}
    \end{figure}

    In manifold optimization, ill-conditioning means that the positive curvatures near a minimizer vary greatly across tangent directions.
    Their largest-to-smallest ratio is the positive Hessian condition number, which enters local convergence estimates for gradient methods~\citep[Thm.~4.20]{boumal2023introduction}.
    A large ratio weakens these estimates and explains the slow zigzagging progress illustrated in Figure~\ref{fig:landscape-routes}.
    The standard remedy changes the Riemannian metric, which determines tangent-vector lengths and the gradient and Hessian spectrum~\citep{absil2008optimization}.
    Existing constructions use Lagrangian Hessians~\citep{mishra2016riemannian}, preconditioned geometry and randomization for generalized orthogonality constraints~\citep{shustin2023riemannian,shustin2024faster}, and approximations to Hessian blocks on product manifolds, including left--right preconditioning~\citep{gao2025optimization}.
    These approaches select a preconditioned metric (\textsf{PM}) on a fixed feasible set, and the resulting spectrum is dictated by that metric choice.

    A complementary question is whether the conditioning can be designed rather than inherited. 
    For linear systems, optimal matrix scaling chooses a scaling matrix to minimize the condition number and admits convex~\citep{braatz1994minimizing} and semidefinite programming~\citep{boyd1994linear,qu2025optimal,jambulapati2023structured} formulations. 
    \citet{dougan2025optimal} establish geodesic convexity of the log-condition number under structured group actions. 
    We transfer this design idea to the constrained problems in this paper. 
    The Hessian spectrum can be controlled through the formulation itself, but an exact design needs an explicit connection between the local Hessian spectrum and the constraint weights, which existing preconditioned metrics do not establish. 
    We derive this connection and obtain an exact linear program.

    Our approach reshapes the constraints rather than the metric. 
    We construct reformulations of the original problem, called preconditioned constraints (\textsf{PC}), each equipped with a prescribed metric. 
    The invertible congruence map of \citet[Prop.~3]{peng2026generalized} transforms~\eqref{eq:quadratic-problem} into
    \begin{equation}\label{eq:GQ-problem}
        \min_{X\in\mathrm{Gq}(P,Q)} f(X)=(h\circ\theta)(X),\qquad
        \mathrm{Gq}(P,Q)=\{X:X^\top QX=P\},
    \end{equation}
    where $P,Q$ are SPD or nonsingular skew-symmetric matrices, and $\theta(X)=\breve QX\breve P^{-1}$ for congruence factors $\breve Q,\breve P$ satisfying~\eqref{eq:congruence}. 
    This map preserves objective values and gives a one-to-one correspondence of minimizers and critical points.
    Each transformed feasible set carries the canonical-like metric of that work.
    Within the resulting family, we minimize the Hessian condition number at a minimizer.

    The matrices $P$ and $Q$ parametrize the family of reshaped constraints.
    Their positive canonical values are the tunable design variables, called the \emph{constraint weights}.
    These spectra consist of the positive eigenvalues of $P$ and $Q$ in the symmetric case and the positive canonical moduli of their $2\times2$ blocks in the skew-symmetric case. 
    Pairing the canonical bases ensures that the fixed spectra of $C$ and $D$ determine the minimizing assignment.
    Figure~\ref{fig:landscape-routes} compares the local geometry produced by a metric change and by constraint reshaping in the Stiefel case.

    In this paper, our contributions are threefold.

    \noindent(1) We develop a paired congruence construction that reshapes the constraints while preserving the underlying quadratic optimization problem.
    We derive the exact positive Hessian spectrum at a minimizer and express its condition number through the constraint weights.
    The spectral decompositions also identify minimizing eigenspaces.
    This construction also admits an indirect interpretation as metric design. Transporting the prescribed metric to the original feasible set yields an explicitly characterized local family.

    \noindent(2) We reduce condition-number minimization to a linear program in the logarithms of the weights and prove that it attains its minimum.
    Its solution gives globally optimal weights within the paired equal-block family.
    LP strong duality gives matching lower bounds on the condition number.

    \noindent(3) We test the reshaped formulations on eight quadratic instances.
    On the orthogonality-constrained problems, we compare them with two preconditioned metric constructions.
    The optimized formulation markedly reduces the local condition number.

    Sections~\ref{sec:preliminaries} and~\ref{sec:critical-geometry} develop the equivalent feasible sets, minimizing orbit, and positive spectrum.
    Section~\ref{sec:condition-number} derives the design LP and dual bounds.
    Sections~\ref{sec:metric-comparison} and~\ref{sec:numerical} present implementation, metric comparisons, and experiments.
    Proofs and detailed metric comparisons appear in the \prooflocation.

\section{Equivalent Reshaped Constraints and Local Curvature}
\label{sec:preliminaries}
    An invertible change of variables preserves the optimization problem.
    The transformed constraints and prescribed metric determine the local geometry.
    We derive the critical-point Hessian and separate it into blocks for motions within the selected column space and into its complement.

    We denote by $\otimes$, $\mathrm{vec}(\cdot)$, $\mathrm{diag}(\cdot)$, and $\mathrm{blkdiag}(\cdot)$ the Kronecker product, columnwise vectorization, diagonal matrix construction, and block-diagonal matrix construction, respectively.
    We use $\langle U,V\rangle=\mathrm{tr}(U^\top V)$ and $\|U\|_{\mathrm{F}}^2=\langle U,U\rangle$, and write $e_i$ for a coordinate vector of the indicated dimension.
    The multiset $\lambda(M)$ contains the ordinary eigenvalues of $M$, including multiplicities, and $\lambda_+(M)$ contains its positive real eigenvalues. We denote positive skew canonical moduli by $\varLambda_M$.
    The groups $\mathrm{GL}(m)$, $\mathrm O(m)$, and $\mathrm U(m)$ comprise invertible real, orthogonal, and unitary $m\times m$ matrices, respectively.

    Throughout, $P\in\mathbb R^{p\times p}$ and $Q\in\mathbb R^{n\times n}$ are nonsingular, with $P^\top=\varrho P$, $Q^\top=\varrho Q$, and $\varrho\in\{+1,-1\}$. The dimensions $p,n$ are even in the skew case.
    Following~\citet[Sec.~3.1]{peng2026generalized}, define $A^{(P,Q)}=P^{-1}A^\top Q$ and $A^{(P)}=P^{-1}A^\top P$.
    The $P$-symmetric and $P$-skew-symmetric subspaces are $\mathcal S_{\mathrm{sym}}(P)=\{A:A^{(P)}=A\}$ and $\mathcal S_{\mathrm{skew}}(P)=\{A:A^{(P)}=-A\}$, with complementary projectors $\mathcal P_{\mathrm{sym}}^{(P)}(A)=(A+A^{(P)})/2$ and $\mathcal P_{\mathrm{skew}}^{(P)}(A)=(A-A^{(P)})/2$.

\subsection{Equivalent Reshaped Constraints}
\label{sec:congruent}
    We restrict to $P,Q\succ0$ when $\varrho=+1$ and nonsingular skew-symmetric forms when $\varrho=-1$. Write
    \begin{equation}\label{eq:spectral-sigma}
    Q=\varPhi_Q\varSigma_Q\varPhi_Q^\top,\qquad
    \varSigma_Q=\begin{cases}\varLambda_Q,&\varrho=+1,\\\varLambda_Q\otimes J_2,&\varrho=-1,\end{cases}
    \end{equation}
    where $\varPhi_Q$ is orthogonal and $\varLambda_Q$ is positive diagonal. Define the analogous quantities for $P$, and set
    \begin{equation}\label{eq:K-rho}
    K_\varrho(Q)=\begin{cases}I_n,&\varrho=+1,\\J_n,&\varrho=-1,\end{cases}
    \qquad K_\varrho(P)=\begin{cases}I_p,&\varrho=+1,\\J_p,&\varrho=-1.\end{cases}
    \end{equation}

    \begin{proposition}\label{prop:diffeo}
    Let $1\le p\le n$, with $P,Q\succ0$ for $\varrho=+1$ or nonsingular skew-symmetric $P,Q$ and even $p,n$ for $\varrho=-1$.
    Let $K_\varrho$ be defined by~\eqref{eq:K-rho}. For invertible factors satisfying
    \begin{equation}\label{eq:congruence}
    Q=\breve Q^\top K_\varrho(Q)\breve Q,\qquad P=\breve P^\top K_\varrho(P)\breve P,
    \end{equation}
    the map $\theta(X)=\breve QX\breve P^{-1}$ is a diffeomorphism from $\mathrm{Gq}(P,Q)$ to $\mathrm{St}(p,n)$ when $\varrho=+1$ and to $\mathrm{Sp}(p,n)$ when $\varrho=-1$, with inverse $\theta^{-1}(Y)=\breve Q^{-1}Y\breve P$.
    \end{proposition}
    Proposition~\ref{prop:diffeo} states the congruence construction of \citet[Prop.~3]{peng2026generalized} using the block-diagonal $J_n$.
    Substitution into the constraint proves the assertion for every pair satisfying~\eqref{eq:congruence}.
    The paired construction also controls the transformed cost.

\subsection{Riemannian Hessian at a Critical Point}
\label{sec:rhess}

    By Proposition~\ref{prop:diffeo}, Problem~\eqref{eq:GQ-problem} takes the form
    \begin{equation}\label{eq:recast}\tag{Gq}
        \min_{X\in\mathrm{Gq}(P,Q)} f(X)=\mathrm{tr}(X^\top AXB),
    \end{equation}
    where
    \begin{equation}\label{eq:equivalent-replacement}
        A=\breve{Q}^{\top}C \breve{Q}, \quad B=\breve{P}^{-1}D\breve{P}^{-\top},
    \end{equation}
    and $(\breve{Q},\breve{P})$ satisfies the congruence conditions~\eqref{eq:congruence}.
    We use the normalized complement and canonical-like geometry of \citet[Secs.~3.4 and 4]{peng2026generalized}. Let $X_\perp$ be a basis for $Q^{-1}\ker(X^\top)$ such that $X_\perp^\top Q^\top Q X_\perp = I_{n-p}$.
    Define $E := [X, X_\perp] \in \mathbb{R}^{n \times n}$ and $O := X_\perp^\top Q X_\perp$.
    The canonical-like metric on $\mathrm{Gq}(P,Q)$ is
    \begin{equation}\label{eq:metric-canonical}
        g_X^{\beta}(Z_1,Z_2)
        = \mathrm{tr}(Z_1^{\top} M_X^{\beta} Z_2)
        \quad\text{with}\quad
        M_X^{\beta} = E^{-\top}\,\mathrm{blkdiag}(\beta PP^\top, OO^\top)\,E^{-1},
        \quad \beta>0.
    \end{equation}
    The ambient extension $\bar f(X) = \mathrm{tr}(X^\top A X B)$ satisfies $\nabla\bar f(X) = 2 A X B$ and $\nabla^2\bar f(X)[Z] = 2 A Z B$.
    Under~\eqref{eq:metric-canonical}, the Riemannian gradient is
    \begin{equation*}
        \mathrm{grad}\,f(X)
        = \mathcal{P}_{\mathrm{T}_X}\bigl((M_X^{\beta})^{-1} \nabla\bar f(X)\bigr),
    \end{equation*}
    where $\mathcal P_{\mathrm{T}_X}(W) = W - X\mathcal{P}_{\mathrm{sym}}^{(P)}(X^{(P, Q)}W)$ is the $g_X^\beta$-orthogonal projection onto the tangent space $\mathrm{T}_X\mathrm{Gq}(P,Q)=\{X\varOmega+X_\perp W:\varOmega\in\mathcal{S}_{\mathrm{skew}}(P),\ W\in\mathbb R^{(n-p)\times p}\}$.
    The set of critical points of \eqref{eq:GQ-problem} is $\mathcal{X} := \{X\in\mathrm{Gq}(P, Q)\,:\, \mathrm{grad}\,f(X) =0\}$.
    The following lemma separates the selected spectrum from its complement at a critical point.
    \begin{lemma}\label{lem:critical-point}
        Let $A,B$ be given by~\eqref{eq:equivalent-replacement} for SPD data $C,D$ and factors~\eqref{eq:congruence}.
        If $X$ is a critical point of~\eqref{eq:recast}, with $X_\perp,O$ as in~\eqref{eq:metric-canonical}, then
        \begin{equation}\label{eq:critical-point-eq}
            X^\top AX BP=PBX^\top AX.
        \end{equation}
        The eigenvalues of $P^{-1}X^\top AX$ form a $p$-element submultiset of $\lambda(Q^{-1}A)$.
        Its complement consists of the eigenvalues of $O^{-1}X_\perp^\top AX_\perp$, equivalently the generalized eigenvalues of the pencil $\bigl(X_\perp^\top AX_\perp,\,X_\perp^\top QX_\perp\bigr)$.
    \end{lemma}

    At $X\in\mathcal X$, stationarity gives $AXB=QXP^{-1}X^\top AXB$.
    With $L_X=P^{-1}X^\top AXB$, condition~\eqref{eq:critical-point-eq} gives $L_X^\top=\varrho L_X$.
    The stationary Lagrangian $f(X)-\mathrm{tr}(L_X(X^\top QX-P))$ has second variation $2\mathrm{tr}(Z^\top AZB-Z^\top QZL_X)$.
    The metric~\eqref{eq:metric-canonical} converts this second variation into the Hessian operator
    \begin{equation*}
        \mathrm{hess}\,f(X)[Z] = 2\mathcal P_{\mathrm{T}_X}((M_X^\beta)^{-1}(AZB - QZP^{-1}X^\top AXB)).
    \end{equation*}
    The Hessian Rayleigh quotient on nonzero tangent vectors is
    \begin{equation}\label{eq:Rayleigh-quotient}
    q(Z)=\frac{2\mathrm{tr}(Z^\top AZB-Z^\top QZP^{-1}X^\top AXB)}{\mathrm{tr}(Z^\top M_X^\beta Z)}.
    \end{equation}

\subsection{Internal and External Hessian Blocks}
\label{sec:Hessian-general}

    We call tangent motions within $\mathrm{range}(X)$ internal and motions into its $Q$-orthogonal complement external.
    Normalizing these coordinates with respect to the metric separates the Hessian eigenvalue problem into two blocks.
    Write $\widehat\Omega\in\mathbb R^{p\times p}$ and $V\in\mathbb R^{(n-p)\times p}$ in the normalized tangent coordinates
    \begin{equation}\label{eq:normalized-tangent}
        Z=\beta^{-1/2}XP^{-1}\widehat\Omega+X_\perp O^{-1}V,
        \qquad \widehat\Omega^\top=-\varrho \widehat\Omega,
        \qquad g_X^\beta(Z,Z)=\|\widehat\Omega\|_{\mathrm{F}}^2+\|V\|_{\mathrm{F}}^2.
    \end{equation}
    Indeed, for $Z=X\varOmega+X_\perp W$, set $\widehat\Omega=\sqrt\beta P\varOmega$ and $V=OW$.
    The identity $\widehat\Omega^\top=-\varrho \widehat\Omega$ is equivalent to $\varOmega^{(P)}=-\varOmega$, and~\eqref{eq:metric-canonical} gives the displayed norm identity.
    Let $T\in\mathbb R^{p^2\times p^2}$ be the commutation matrix defined by $T\mathrm{vec}(\widehat\Omega)=\mathrm{vec}(\widehat\Omega^\top)$ for every $\widehat\Omega$. It satisfies $T^\top=T$ and $T^2=I$. The orthogonal projector onto the internal coordinate space is $\Pi_\varrho=(I-\varrho T)/2$, and this space has dimension $p(p-\varrho)/2$.

    \begin{proposition}\label{prop:Rayleigh-quotient}
        For the data and factors in Lemma~\ref{lem:critical-point}, let $X$ be a critical point of~\eqref{eq:recast} with metric~\eqref{eq:metric-canonical}, $\beta>0$, and normalized coordinates~\eqref{eq:normalized-tangent}.
        Let $T$ be the commutation matrix on $p\times p$ coordinates.
        The nonzero stationary values of~\eqref{eq:Rayleigh-quotient} are the nonzero intrinsic Hessian eigenvalues on $\mathrm{T}_X\mathrm{Gq}(P,Q)$.
        They coincide with the nonzero eigenvalues of the real $np\times np$ reduced-coordinate matrix
        \begin{equation}\label{eq:Hf}
            H_f := \mathrm{blkdiag}\bigl(\frac{1}{\beta}(I_{p^2} - \varrho T)\,H_{\mathrm{int}},\; 2H_{\mathrm{ext}}\bigr),
        \end{equation}
        with the internal block and external block given by
        \begin{equation*}
            \begin{aligned}
                H_{\mathrm{int}} &:= B\otimes P^{-\top} X^{\top} A X P^{-1} - P^{-1} X^{\top} A X B\otimes P^{-1},\\
                H_{\mathrm{ext}} &:= B\otimes O^{-\top} X_\perp^{\top} A X_\perp O^{-1} - P^{-1} X^{\top} A X B\otimes O^{-1}.
            \end{aligned}
        \end{equation*}
        Restricting the internal block to $\mathrm{vec}\{U:U^\top=-\varrho U\}$ gives the intrinsic spectrum. Internal coordinates outside this subspace contribute only redundant zeros.
    \end{proposition}

    \begin{remark}\label{rem:Hf-decoupling}
        The parameter $\beta$ in Proposition~\ref{prop:Rayleigh-quotient} scales only the internal Hessian block. The external block is independent of $\beta$.
    \end{remark}

\section{Realizable Constraint Weights and Positive Curvature}
\label{sec:critical-geometry}

    For fixed data $C,D$, we construct constraint matrices $P,Q$ and congruence factors that determine $A,B$.
    Pairing their canonical bases preserves the minimizing assignment and gives explicit formulas for the positive Hessian eigenvalues in terms of the design weights.

\subsection{The Paired Construction}
\label{sec:assumption}
    We require the transformed cost and constraint matrices to satisfy the \emph{commutation conditions}
    \begin{equation}\label{eq:commutation-conditions}
        BP = PB \qquad\text{and}\qquad AQ = QA,
    \end{equation}
    which allow common canonical bases for $A,Q$ and for $B,P$.
    To express the paired construction in terms of the data $C,D$, define
    \begin{equation}\label{eq:frak-CD}
        \mathfrak{C} := C^{-1/2}\,K_\varrho(Q)\,C^{-1/2},
        \qquad
        \mathfrak{D} := D^{1/2}\,K_\varrho(P)\,D^{1/2}.
    \end{equation}
    We use pairs $(P,Q)$ with \emph{positive canonical spectra}, defined by $[\varLambda_P]_i>0$ and $[\varLambda_Q]_j>0$.
    These conditions require $P,Q\succ0$ for $\varrho=+1$ and positive skew canonical moduli for $\varrho=-1$.
    In the symmetric case, $K_\varrho(P)=I_p$ and the pair $(\mathfrak C,\mathfrak D)=(C^{-1},D)$ is fixed.
    In the skew-symmetric case $K_\varrho(P)=J_p$, $\mathfrak C,\mathfrak D$ are skew-symmetric, and $[\varLambda_{\mathfrak C}]_i,[\varLambda_{\mathfrak D}]_i$ denote the positive canonical moduli of length $n/2$ and $p/2$ collected in the form~\eqref{eq:spectral-sigma}.
    The numbers of canonical moduli of $P$ and $Q$ are
    \begin{equation}\label{eq:canonical-dimensions}
        (\tilde p,\tilde n)=
        \begin{cases}
            (p,n), & \varrho=+1,\\
            (p/2,n/2), & \varrho=-1.
        \end{cases}
    \end{equation}
    Write $c_i=[\varLambda_{\mathfrak C}]_i$ and $d_i=[\varLambda_{\mathfrak D}]_i$.
    Fix the data ordering
    \begin{equation}\label{eq:ordering}
        c_1\ge\cdots\ge c_{\tilde n}>0,\qquad d_1\ge\cdots\ge d_{\tilde p}>0.
    \end{equation}
    The positive design moduli $[\varLambda_P]_i,[\varLambda_Q]_j$ retain these paired labels and need not be monotone.
    Canonical matrices use $J_2$ blocks in the skew case. For such matrices, $[\cdot]_{\mathcal I}$ selects the ordered blocks in $\mathcal I$, whereas $[\cdot]_{\le p}$ selects the leading $p$ ordinary coordinates.

\begin{remark}\label{rem:genericity}
Throughout, we assume $c_1>\cdots>c_{\tilde n}>0$ in~\eqref{eq:ordering}. Repeated $d_i$ are allowed. The equal-block design restriction is formalized in Corollary~\ref{cor:diagonal-block}. We omit empty index families and require at least one positive Hessian mode.
\end{remark}
At a critical point, the blocks of Proposition~\ref{prop:Rayleigh-quotient} also satisfy
\begin{equation*}
\begin{aligned}
H_{\mathrm{int}}&=(I_p\otimes P^{-\top}X^\top AX-P^{-1}X^\top AX\otimes I_p)(B\otimes P^{-1}),\\
H_{\mathrm{ext}}&=(I_p\otimes O^{-\top}X_\perp^\top AX_\perp-P^{-1}X^\top AX\otimes I_{n-p})(B\otimes O^{-1}).
\end{aligned}
\end{equation*}

    \begin{proposition}\label{prop:commutation}
        Fix SPD data $C,D$ and $\varrho\in\{+1,-1\}$, with even $p,n$ for $\varrho=-1$.
        Let $\varSigma_{\mathfrak C},\varSigma_{\mathfrak D}$ be the canonical forms of~\eqref{eq:frak-CD} in the order~\eqref{eq:ordering}.
        For any positive design moduli in $\varSigma_P,\varSigma_Q$ and orthogonal canonical bases $\varPhi_P,\varPhi_Q$, there exist congruence factors $\breve{Q},\breve{P}$ satisfying~\eqref{eq:congruence}.
        The resulting $A,B$ in~\eqref{eq:equivalent-replacement} have the following properties.
        \begin{enumerate}
            \item $A,Q,A^{-1}Q$ share the basis $\varPhi_Q$, and $B,P,BP$ share the basis $\varPhi_P$. In particular,~\eqref{eq:commutation-conditions} holds.
            \item $\lambda(A^{-1}Q) = \lambda(\mathfrak{C})$ and $\lambda(BP) = \lambda(\mathfrak{D})$.
            \item $\varSigma_A = \varSigma_Q\varSigma_{\mathfrak{C}}^{-1}$ and $\varSigma_B = \varSigma_P^{-1}\varSigma_{\mathfrak{D}}$.
        \end{enumerate}
    \end{proposition}

\paragraph{Constructing the Factors.}\label{rem:construction}
Decompose $\mathfrak C=\varPhi_{\mathfrak C}\varSigma_{\mathfrak C}\varPhi_{\mathfrak C}^\top$ and $\mathfrak D=\varPhi_{\mathfrak D}\varSigma_{\mathfrak D}\varPhi_{\mathfrak D}^\top$ in the paired order~\eqref{eq:ordering}. The required factors are
\begin{equation}\label{eq:breve-construction}
\begin{aligned}
\breve Q&=C^{-1/2}\varPhi_{\mathfrak C}(\varSigma_Q\varSigma_{\mathfrak C}^{-1})^{1/2}\varPhi_Q^\top,\\
\breve P&=D^{1/2}\varPhi_{\mathfrak D}(\varSigma_P\varSigma_{\mathfrak D}^{-1})^{1/2}\varPhi_P^\top.
\end{aligned}
\end{equation}
Both products under the square roots are SPD. Thus, for arbitrary SPD data $C,D$, the factors realize every positive design and every orthogonal basis pair $\varPhi_P,\varPhi_Q$.

\subsection{The Minimizing Orbit}\label{sec:minimizing-orbit}
For a canonical form $\varSigma_P$, write $\varSigma_P^{1/2}=\varLambda_P^{1/2}$ in the symmetric case and $\varLambda_P^{1/2}\otimes I_2$ in the skew case. This denotes the positive modulus root, so that $\varSigma_P=\varSigma_P^{1/2}K_\varrho(P)\varSigma_P^{1/2}$. Likewise for $Q$.
\begin{proposition}\label{prop:critical-point}
For the paired factors of Proposition~\ref{prop:commutation} and the ordered data of Remark~\ref{rem:genericity}, an aligned minimizer of~\eqref{eq:recast} is
\begin{equation}\label{eq:critical-point}
X_\star=[\varPhi_Q]_{:,\le p}[\varSigma_Q]_{\le p}^{-1/2}\varSigma_P^{1/2}\varPhi_P^\top.
\end{equation}
Its ordered Gram matrix is $\varPhi_P^\top X_\star^\top AX_\star\varPhi_P=\varSigma_P[\varSigma_{\mathfrak C}]_{\le p}^{-1}$. The minimum is independent of the design and satisfies
\begin{equation*}
f_\star=h_\star=\begin{cases}\sum_{i=1}^{p}d_i/c_i,&\varrho=+1,\\2\sum_{i=1}^{p/2}d_i/c_i,&\varrho=-1.\end{cases}
\end{equation*}
\end{proposition}
The minimizing assignment and value are independent of the design weights.
Symmetries of the minimizer set determine the zero-curvature directions that the positive condition number excludes.
Define
\begin{equation}\label{eq:gauge-group}
\mathcal U=\{U\in\mathrm{GL}(p):U^\top PU=P,\ UBU^\top=B\}.
\end{equation}

\begin{proposition}\label{prop:critical-manifold}
Use the paired construction of Proposition~\ref{prop:commutation}, with $c_1>\cdots>c_{\tilde n}>0$ in~\eqref{eq:ordering}.
Let $\delta_1>\cdots>\delta_s>0$ be the distinct values among $d_1,\ldots,d_{\tilde p}$.
Their blocks $\mathcal D_b=\{i:d_i=\delta_b\}$ have sizes $m_b=|\mathcal D_b|$, with $\sum_{b=1}^s m_b=\tilde p$ as in~\eqref{eq:canonical-dimensions}.
The minimizer set of~\eqref{eq:recast} is the compact embedded orbit $\mathcal O=X_\star\mathcal U$, with $X_\star$ from~\eqref{eq:critical-point} and $\mathcal U$ from~\eqref{eq:gauge-group}.
For every $X\in\mathcal O$, $\mathcal O=X\mathcal U$.
The group $\mathcal U$ is isomorphic to $\prod_{b=1}^s\mathrm O(m_b)$ for $\varrho=+1$ and $\prod_{b=1}^s\mathrm U(m_b)$ for $\varrho=-1$.
Under~\eqref{eq:metric-canonical}, the intrinsic Hessian satisfies
\[
 \ker\bigl(\operatorname{hess}f(X)\bigr)=\mathrm T_X\mathcal O,
 \qquad X\in\mathcal O.
\]
Its kernel has dimension
\begin{equation}\label{eq:dim-centralizer}
\dim\mathcal U=\begin{cases}\sum_{b=1}^s m_b(m_b-1)/2,&\varrho=+1,\\\sum_{b=1}^s m_b^2,&\varrho=-1.\end{cases}
\end{equation}
The Hessian is positive definite on the $g_X^\beta$-normal complement of $\mathrm T_X\mathcal O$ in $\mathrm T_X\mathrm{Gq}(P,Q)$.
\end{proposition}
When the $d_i$ are distinct, the symmetric minimizing orbit consists of $2^p$ isolated sign choices.
In the skew case, phase rotations of the selected two-planes form the torus $\mathrm U(1)^{p/2}$ of dimension $p/2$.
Repeated $d_i$ also permit rotations within equal-data blocks.
\begin{corollary}\label{cor:diagonal-block}
Under the hypotheses of Proposition~\ref{prop:critical-manifold}, suppose $[\varLambda_P]_i=[\varLambda_P]_{i'}$ whenever $d_i=d_{i'}$.
The action of $\mathcal U$ in~\eqref{eq:gauge-group} is free, proper, and isometric for~\eqref{eq:metric-canonical}. Thus $\mathrm{Gq}(P,Q)/\mathcal U$ is a Riemannian quotient with unique minimizer $[X_\star]$. The tangent at $[X_\star]$ identifies with the metric-normal complement of $\mathrm T_X\mathcal O$ for $X\in\mathcal O$, and its Hessian spectrum is each representative's positive spectrum.
\end{corollary}
The equal-block restriction makes the condition number constant along the minimizing orbit and imposes linear equalities on the logarithmic design variables.
The distinct values $c_i$ fix the selected invariant subspaces. Repeated values require additional subspace choices.

\subsection{Positive Curvature Modes}
\label{sec:spectrum}

We compare all induced metrics on a fixed reference quotient.
Put $K_n=K_\varrho(Q)$, $K_p=K_\varrho(P)$, and $\nu=1$ or $2$ in the symmetric or skew case, respectively.
With the positive modulus roots defined above, set
\begin{equation}\label{eq:common-reference}
 W=\varSigma_Q^{1/2}\varPhi_Q^\top X\varPhi_P\varSigma_P^{-1/2},\qquad
 C_0=\mathrm{diag}(c_i^{-1})\otimes I_\nu,\quad
 D_0=\mathrm{diag}(d_i)\otimes I_\nu.
\end{equation}
For every design, $W^\top K_nW=K_p$ and $f(X)=\mathrm{tr}(W^\top C_0WD_0)$.
The map is equivariant under the fixed compact group
$\mathcal C_D=\{R\in\mathrm O(p):R^\top K_pR=K_p,\ RD_0=D_0R\}$.
This equivariance identifies all realization quotients with
$\mathcal Q=\{W:W^\top K_nW=K_p\}/\mathcal C_D$.
Their common minimizing point is $[W_\star]$, where $W_\star=[I_p,0]^\top$.
Internal motions between distinct $d$ blocks remain in the quotient. For scalar $D_0$, it is the Grassmann quotient in the symmetric case and a quotient by $\mathrm U(p/2)$ in the skew case.

\begin{theorem}\label{thm:eigenvalue-Hf}
Under Proposition~\ref{prop:commutation}, Remark~\ref{rem:genericity}, and the equal-block restriction of Corollary~\ref{cor:diagonal-block}, use metric~\eqref{eq:metric-canonical} with $\beta>0$.
Set $\alpha_i=[\varLambda_P]_i$ and $\gamma_j=[\varLambda_Q]_j$.
In one design-independent basis of $\mathrm T_{[W_\star]}\mathcal Q$, the second differential and the induced metric have SPD matrices
\begin{equation}\label{eq:quotient-SPD}
 H_0=\mathrm{diag}(b_k),\qquad G=\mathrm{diag}(g_k),\qquad
 \mathrm{hess}\,\bar f([W_\star])=G^{-1}H_0,
\end{equation}
where $\bar f$ is the objective on $\mathcal Q$.
Use $k=(i,j,\mathrm{int},s)$ for an internal mode and $k=(i,j,\mathrm{ext},s)$ for an external mode, where $s\in\{-1,+1\}$. Use $k=(i,0)$ for a diagonal mode.
The entries in~\eqref{eq:quotient-SPD} are
\begin{equation}\label{eq:eig-symm}
 (b_k,g_k)=
 \begin{cases}
 \left(\dfrac{(d_i+s d_j)(c_i+s c_j)}{c_ic_j},\,\beta\alpha_i\alpha_j\right),
 &k=(i,j,\mathrm{int},s),\quad 1\le i<j\le\tilde p,\\[6pt]
 \left(\dfrac{2d_i(c_i+s c_j)}{c_ic_j},\,\dfrac{\alpha_i}{\gamma_j}\right),
 &k=(i,j,\mathrm{ext},s),\quad 1\le i\le\tilde p<j\le\tilde n,\\[6pt]
 \left(\dfrac{4d_i}{c_i},\,\beta\alpha_i^2\right),
 &k=(i,0),\quad 1\le i\le\tilde p.
 \end{cases}
\end{equation}
Retain $s=-1$ modes, excluding internal modes with $d_i=d_j$. The skew case also includes all $s=+1$ and diagonal modes. Each label contributes one entry per matrix in~\eqref{eq:quotient-SPD} in the symmetric case and two identical entries in the skew case.
\end{theorem}
\begin{remark}\label{rem:spectrum-consequence}
The metric inverse $G^{-1}$ converts the fixed curvature form $H_0$ into the Hessian operator, and the symmetric SPD matrix $G^{-1/2}H_0G^{-1/2}$ has the same spectrum. The ratios $\lambda_k=b_k/g_k$ therefore give the complete positive intrinsic Hessian spectrum on the minimizing orbit, counting multiplicities even when values coincide. Neither the leading weights $\gamma_j$ with $j\le\tilde p$ nor the orthogonal realization bases affect the pair of forms.

Varying positive equal-block weights exhausts the paired metric family at $[W_\star]$. Without the equal-block restriction, the ratios give the positive spectrum at the aligned minimizer.
\end{remark}

\paragraph{Condition Number.}\label{rem:zero-eigs}
For the nonempty multiset $\mathcal S=\{b_k/g_k\}$ of retained modes, the common positive condition number is
\begin{equation*}
 \kappa_{\mathrm{align}}=\kappa(G^{-1}H_0)
 =\frac{\max\mathcal S}{\min\mathcal S}.
\end{equation*}
The quotient removes the intrinsic zero modes, which span the orbit tangent.
The extended internal matrix~\eqref{eq:Hf} has a further $p(p+\varrho)/2$ redundant zeros outside the intrinsic Hessian domain.

    \begin{remark}\label{rem:rescaling}
        For fixed $C,D,\beta$ under Theorem~\ref{thm:eigenvalue-Hf}, the aligned spectrum $\lambda(H_f)$ depends only on $(\varLambda_P, \varLambda_Q)$ and not on $(\varPhi_P, \varPhi_Q)$.
        Every eigenvalue in Theorem~\ref{thm:eigenvalue-Hf} scales by the common factor $t/s$ under the joint rescaling $\varLambda_P\mapsto s \, \varLambda_P$, $\varLambda_Q\mapsto t\,\varLambda_Q$, $\beta\mapsto \beta\,(st)^{-1}$ for $s,t>0$.
        The invariance of $\kappa(\mathrm{hess}\,f)$ under this rescaling gives the triple $(\varLambda_P,\varLambda_Q,\beta)$ a two-parameter redundancy.
        In particular, we may fix $\beta=1$.
    \end{remark}

    \begin{remark}\label{rem:special-cases}
        Under Theorem~\ref{thm:eigenvalue-Hf} with $\beta=1$, the \emph{identity design} sets every constraint modulus to one, $\varLambda_P=\varLambda_Q=I$.
        Theorem~\ref{thm:eigenvalue-Hf} gives the aligned Stiefel-type spectrum determined by $(\mathfrak C,\mathfrak D)$.
        The \emph{balanced heuristic} sets the constraint moduli equal to the corresponding data moduli, $\varLambda_P=\varLambda_{\mathfrak D}$ and $\varLambda_Q=\varLambda_{\mathfrak C}$, making $A=B=I$ under~\eqref{eq:breve-construction}.
        Its external eigenvalues are $2(1-c_j/c_i)$ in the symmetric case and $2(1\pm c_j/c_i)$ in the skew case, lying in $[0,2]$ and $[0,4]$, respectively.
        For $p<n$, these bounds control the largest external eigenvalue, while the smallest positive value can approach zero as the data vary.
        The internal spectral spread also depends on the data.
    \end{remark}

\section{Minimizing the Positive Hessian Condition Number}
\label{sec:condition-number}

    Each positive eigenvalue in Theorem~\ref{thm:eigenvalue-Hf} is a fixed data coefficient multiplied by products and ratios of positive design weights.
    Its logarithm is affine in the logarithmic weights, as for monomials in geometric programming~\citep{boyd2007tutorial}.
    The width of the smallest interval containing these logarithms equals the logarithm of the positive Hessian condition number.
    We minimize this width by an LP and derive matching lower bounds from its dual.

\subsection{The Design Linear Program}
\label{sec:design-LP}

    With $c,d$ fixed in the order~\eqref{eq:ordering}, introduce the logarithmic design variables
    \begin{equation*}        x_i := \log[\varLambda_P]_i \ (1\le i\le \tilde p),
        \qquad
        y_j := \log[\varLambda_Q]_j \ (\tilde p< j\le \tilde n),
        \qquad
        \tau := \log\beta,
    \end{equation*}
    with the canonical dimensions $(\tilde p,\tilde n)$ defined in~\eqref{eq:canonical-dimensions}. The leading $\tilde p$ moduli of $Q$ do not enter the aligned spectrum and may take any positive values.
    Exclude the zero internal minus modes with $d_i=d_j$ before taking logarithms.
    The offsets $a^{\pm}_{ij}$, $b^{\pm}_{ij}$, and $a^0_i$ are the logarithms of the retained curvature coefficients $b_k$ in~\eqref{eq:eig-symm}. Define
    \begin{equation*}
        a^{\pm}_{ij} := \log\frac{(d_i\pm d_j)(c_i\pm c_j)}{c_ic_j}, \qquad
        b^{\pm}_{ij} := \log\frac{2\,d_i\,(c_i\pm c_j)}{c_ic_j}, \qquad
        a^0_i :=\log(4d_i/c_i).
    \end{equation*}
    The variables $\overline\lambda,\underline\lambda$ bound the log-eigenvalues from above and below.

    \textbf{Symmetric case ($\varrho=+1$)}.
    The log-eigenvalues are
    \begin{equation*}
        \ell^{\nu}_{ij}=a^-_{ij}-x_i-x_j-\tau\ \ (i<j\le \tilde p),
        \qquad
        \ell^{\mu}_{ij}=b^-_{ij}+y_j-x_i\ \ (i\le \tilde p<j\le \tilde n).
    \end{equation*}
    \textbf{Skew-symmetric case ($\varrho=-1$)}.
    The log-eigenvalues are
    \begin{equation*}
        \begin{aligned}
            \ell^{\nu,\pm}_{ij}&=a^{\pm}_{ij}-x_i-x_j-\tau\ (i<j\le \tilde p), &
            \ell^{0}_{i}&=a^0_i-2x_i-\tau\ (i\le \tilde p),\\
            \ell^{\mu,\pm}_{ij}&=b^{\pm}_{ij}+y_j-x_i\ (i\le \tilde p<j\le \tilde n). &&
        \end{aligned}
    \end{equation*}
Let $\mathcal E(x,y,\tau)$ collect the retained affine rows for the chosen symmetry case.
Both cases give the linear program
\begin{equation}\label{eq:LP-symm}\tag{LP}
 \begin{aligned}
 \operatorname*{minimize}_{x,y,\tau,\overline\lambda,\underline\lambda}
 &\quad \overline\lambda-\underline\lambda\\
 \text{subject to}\quad
 &\underline\lambda\le\ell\le\overline\lambda,
       &&\ell\in\mathcal E(x,y,\tau),\\
 &x_i=x_{i'},&&d_i=d_{i'}.
 \end{aligned}
\end{equation}
\phantomsection\label{sec:worst-case-rationale}
The equal-block constraints make the metric descend to the quotient in Corollary~\ref{cor:diagonal-block}.
The factors~\eqref{eq:breve-construction} realize every feasible design. The optimal value determines
\begin{equation}\label{eq:kappa-star}
 \kappa_\star=\exp(\overline\lambda^\star-\underline\lambda^\star).
\end{equation}
The program has $O(\tilde p+\tilde n)$ variables and $O(\tilde p^2+\tilde p(\tilde n-\tilde p))$ inequalities.

    \begin{remark}\label{rem:degenerate-design}
        We remove the identically zero modes of Theorem~\ref{thm:eigenvalue-Hf} before forming~\eqref{eq:LP-symm}.
        When the $c_i$ are distinct, repeated $d_i$ give zero symmetric internal modes and skew internal minus modes. All these modes are tangent to the minimizing orbit.
        The skew plus and diagonal modes remain positive.
        A nonempty positive spectrum is required to define the condition number.
        The external modes ensure this for $n>p$. When $n=p$, a positive internal mode must remain.
    \end{remark}

\subsection{Well-Posedness and Invariances}
\label{sec:LP-properties}

    We eliminate the equality constraints by introducing one variable $\hat x_b$ for each equal-$d$ block $\mathcal D_b$ in Proposition~\ref{prop:critical-manifold}.
    For example, $d=(2,2,1)$ requires $x_1=x_2$.
    The independent variables and the corresponding full vector are
    \begin{equation*}
        z=(\hat x_1,\hat x_2,y_4,\tau),\qquad
        (x,y,\tau)=(\hat x_1,\hat x_1,\hat x_2,y_4,\tau).
    \end{equation*}
    In general, write $z=(\hat x,y,\tau)\in\mathbb R^r$ and $(x,y,\tau)=Rz$.
    The matrix $R$ repeats $\hat x_b$ at each position $i\in\mathcal D_b$ and leaves $y$ and $\tau$ unchanged.
    This enforces the equal-block constraints for every $z$.

    Count the retained modes with their multiplicities in Theorem~\ref{thm:eigenvalue-Hf}.
    Each log-eigenvalue is an affine function of these independent variables.
    Index these functions by $k\in\mathcal K$ and write
    \begin{equation}\label{eq:affine-log-modes}
        \ell_k(z)=a_k+\langle\xi_k,z\rangle,
        \qquad \xi_k=R^\top\bar\xi_k,
    \end{equation}
    where $a=(a_k)_{k\in\mathcal K}$ collects the constant terms and $\bar\xi_k$ is the coefficient vector before substitution of $(x,y,\tau)=Rz$.
    We use these reduced coefficients throughout the LP analysis.
    Let $\Xi$ have rows $\xi_k^\top$.
    At the quotient minimizer, Theorem~\ref{thm:eigenvalue-Hf} identifies the exact family of realizable metric matrices as
    \begin{equation}\label{eq:exact-metric-family}
      \mathcal G_{\text{\normalfont\sffamily PC}}
      =\{\mathrm{diag}(e^{-(\Xi z)_k}):z\in\mathbb R^r\},
      \qquad H_0=\mathrm{diag}(e^{a_k}).
    \end{equation}
    The LP therefore minimizes the positive condition number over precisely this metric family.
    A positive diagonal matrix $\mathrm{diag}(m_k)$ belongs to the family if and only if
    $\sum_k w_k\log m_k=0$ for every $w\in\ker\Xi^\top$.
    These identities express the restrictions imposed by shared constraint weights.
    The matrix representations retain all mode multiplicities. When solving the LP numerically, identical affine functions, including their offsets, need only be constrained once.

    \begin{lemma}        \label{lem:LP-invariance}
        For the linear program~\eqref{eq:LP-symm}, with zero modes removed as in Remark~\ref{rem:degenerate-design}, the spread is invariant under the joint shift
        \begin{equation*}
            x_i\mapsto x_i+\bar a\ (\forall i), \qquad
            y_j\mapsto y_j+\bar b\ (\forall j), \qquad
            \tau\mapsto \tau+t, \qquad\text{with}\quad \bar b=-\bar a-t,
        \end{equation*}
        which translates every log-eigenvalue by the common constant $-(2\bar a+t)$.
        In particular, we may fix $\beta=1$ and normalize one canonical spectral degree of freedom without loss of optimality.
    \end{lemma}

    \begin{proposition}        \label{prop:LP-wellposed}
        Assume Remark~\ref{rem:genericity} and let $\mathcal K$ index the positive modes retained under Remark~\ref{rem:degenerate-design}.
        Suppose $\mathcal K$ is nonempty.
        Then the linear program~\eqref{eq:LP-symm} is feasible, bounded below by $0$, and attains its minimum.
        A global minimizer $(\varLambda_P^\star,\varLambda_Q^\star)$ therefore exists, and \eqref{eq:kappa-star} defines $\kappa_\star$.
    \end{proposition}

\subsection{Duality and Lower Bounds for the Optimal Spread}
\label{sec:LP-lowerbounds}

    Shared design coefficients can prevent the positive eigenvalues from becoming equal.
    Two affine log-eigenvalues with identical coefficients have a fixed difference, whose absolute value is a lower bound on every feasible interval width.
    The same argument applies to two convex combinations of log-eigenvalues with equal aggregate design coefficients.
    Their difference depends only on the data.
    LP duality shows that the best bound from these cancellations equals the optimal spread.
    Using the reduced affine functions $\ell_k(z)=a_k+\langle\xi_k,z\rangle$, the linear program~\eqref{eq:LP-symm} takes the reduced form
    \begin{equation}\label{eq:primal-spread}\tag{P}
        G_\star=\min_{z,\overline\lambda,\underline\lambda}\ \overline\lambda-\underline\lambda
        \quad\text{s.t.}\quad
        \langle\xi_k,z\rangle-\overline\lambda\le-a_k,\quad
         -\langle\xi_k,z\rangle+\underline\lambda\le a_k,\quad k\in\mathcal{K},
    \end{equation}
    where $z$, $\overline\lambda$, and $\underline\lambda$ are free variables and $\kappa_\star=e^{G_\star}$.

    Multipliers $u_k,v_k\ge0$ for the upper and lower constraints give stationarity conditions $\sum_k u_k=\sum_k v_k=1$ and $\sum_k(u_k-v_k)\xi_k=0$. Thus the dual is
    \begin{equation}\label{eq:dual-spread}\tag{D}
        \max_{u,v\ge0}\ \sum_k(u_k-v_k)a_k
        \quad\text{s.t.}\quad
        \sum_k u_k=1,\quad \sum_k v_k=1,\quad \sum_k(u_k-v_k)\xi_k=0.
    \end{equation}

    Define $w:=u-v$.
    Its positive and negative parts have equal total mass, and the weighted design coefficients cancel.
    Let
    \begin{equation}\label{eq:balanced-weights}
        \mathcal Z:=\Bigl\{w\in\mathbb R^{\mathcal K}:\textstyle\sum_k w_k=0,\ \sum_k w_k\xi_k=0\Bigr\},
        \qquad
        \sigma(w):=\sum_{k:\,w_k>0}w_k=\sum_{k:\,w_k<0}|w_k|.
    \end{equation}
    Equivalently, $\sum_k w_k\bar\xi_k$ is orthogonal to the design subspace $\mathrm{range}(R)$. Every $w\in\mathcal Z$ satisfies $\sum_k w_k\ell_k(z)\equiv\langle w,a\rangle$, independent of $z$. These identities give combinations of log-eigenvalues that are unchanged throughout the realizable metric family~\eqref{eq:exact-metric-family}.

    \begin{theorem}        \label{thm:strong-duality}
        Let $\mathcal K$ be the finite nonempty family in~\eqref{eq:affine-log-modes}, with $\mathcal Z,\sigma$ defined in~\eqref{eq:balanced-weights}.
        The optimal spread in~\eqref{eq:primal-spread} satisfies
        \begin{equation}\label{eq:strong-duality}
            G_\star
            =\max_{w\in\mathcal Z,\ \sigma(w)\le1}\langle w,a\rangle
            =\sup_{w\in\mathcal Z,\ \sigma(w)>0}\frac{\langle w,a\rangle}{\sigma(w)} ,
        \end{equation}
        and both the primal \eqref{eq:primal-spread} and the dual \eqref{eq:dual-spread} are attained.
        The second supremum has an empty domain precisely when $\mathcal Z=\{0\}$.
        In that case, $G_\star=0$ and we use the convention $\sup\emptyset=0$.
    \end{theorem}

    \begin{corollary}        \label{cor:certificate}
        For the affine family and balanced weights of Theorem~\ref{thm:strong-duality}, every $w\in\mathcal Z$ with $\sigma(w)>0$ gives the lower bound $G_\star\ge|\langle w,a\rangle|/\sigma(w)$.
        Equality holds at any nonzero maximizer of~\eqref{eq:strong-duality} normalized by $\sigma(w)=1$, when such a maximizer exists.
    \end{corollary}

    \begin{corollary}        \label{cor:gain}
        Under Theorem~\ref{thm:eigenvalue-Hf}, retain the nonempty positive spectrum and impose the equal-block constraints of Corollary~\ref{cor:diagonal-block}.
        The identity design $(x,y,\tau)=(0,0,0)$ is feasible for~\eqref{eq:LP-symm}.
        It gives $(\varLambda_P,\varLambda_Q)=(I,I)$ and $\beta=1$, the Stiefel-type reformulation of Remark~\ref{rem:special-cases}.
        Hence $\kappa_\star$ does not exceed the aligned condition number of this identity design, with strict inequality unless $(I,I)$ already minimizes the spread.
        For $\varrho=+1$ the conditions~\eqref{eq:commutation-conditions} hold for Problem~\eqref{eq:quadratic-problem} on $\mathrm{St}(p,n)$.
        The identity-design value therefore equals the actual Hessian condition number there, and $\kappa_\star\le\kappa(\mathrm{hess}\,h(Y))$ at the minimizer $Y$.
        For $\varrho=-1$ the conditions~\eqref{eq:commutation-conditions} fail for generic $(C,D)$ on $\mathrm{Sp}(p,n)$.
        The identity-design value belongs to the aligned reformulation on $\mathrm{Sp}(p,n)$ produced by~\eqref{eq:breve-construction}, whose transformed cost differs from the original one. Their Hessian condition numbers need not agree.
    \end{corollary}

\subsection{A Worked Spectral Design}
\label{sec:worked-design}

    This example shows how a repeated data weight produces both an orbit direction and an equal-block design restriction.
    The restriction fixes the difference between two log-eigenvalues and determines the optimal condition number.

    Take $\varrho=+1$, $n=4$, $p=3$, with
\begin{equation*}
C=\mathrm{diag}(1/4,1/2,1,2),\quad D=\mathrm{diag}(2,2,1),\quad c=(4,2,1,1/2),\quad d=(2,2,1).
\end{equation*}
    To identify the orbit direction, let $Y_\star =[e_1,e_2,e_3]$, and let $R_{ij}(t)$ be the identity with the rotation
    $\left[\begin{smallmatrix}\cos t&-\sin t\\\sin t&\cos t\end{smallmatrix}\right]$ in coordinates $i,j$.
    The two column rotations give
    \begin{equation*}
        h(Y_\star R_{12}(t))=\frac52,\qquad
        h(Y_\star R_{13}(t))=\frac52+\frac34\sin^2t.
    \end{equation*}
    The equal data weights keep the first curve on the minimizing orbit. Its tangent spans the one-dimensional Hessian kernel.
    The second curve preserves the selected span while reassigning its eigenvectors across distinct weights.
    At $t=\pi/2$ it has assignment $(3,2,1)$ and objective value $13/4$.

    For an admissible equal-block design, the repeated weight $d_1=d_2$ requires $x_1=x_2$.
    The internal mode $(1,2)$ vanishes and is excluded from the LP.
    Write $\nu_{ij}$ and $\mu_{ij}$ for the internal and external $s=-1$ eigenvalues in~\eqref{eq:eig-symm}, respectively.
    The two positive internal and three external modes give the five affine rows
    \begin{equation}\label{eq:worked-rows}
        \begin{array}{c|c|c}
            k & a_k & \ell_k \\ \hline
            1\ (\nu_{13}) & \log(3/4) & \log(3/4)-x_1-x_3-\tau\\
            2\ (\nu_{23}) & \log(1/2) & \log(1/2)-x_2-x_3-\tau\\
            3\ (\mu_{14}) & \log 7 & \log 7+y_4-x_1\\
            4\ (\mu_{24}) & \log 6 & \log 6+y_4-x_2\\
            5\ (\mu_{34}) & \log 2 & \log 2+y_4-x_3
        \end{array}
    \end{equation}
    The external rows include the factor $2$ in~\eqref{eq:eig-symm}.
    Under $x_1=x_2$, the design coefficients of the first two rows agree, giving the fixed difference
    \begin{equation*}
        \ell_1-\ell_2=\log(3/2),
        \qquad G_\star\ge\log(3/2).
    \end{equation*}
    Thus every design in this family has a condition number of at least $3/2$.

\begin{figure}[htbp]
    \centering
    \includegraphics[width=14cm]{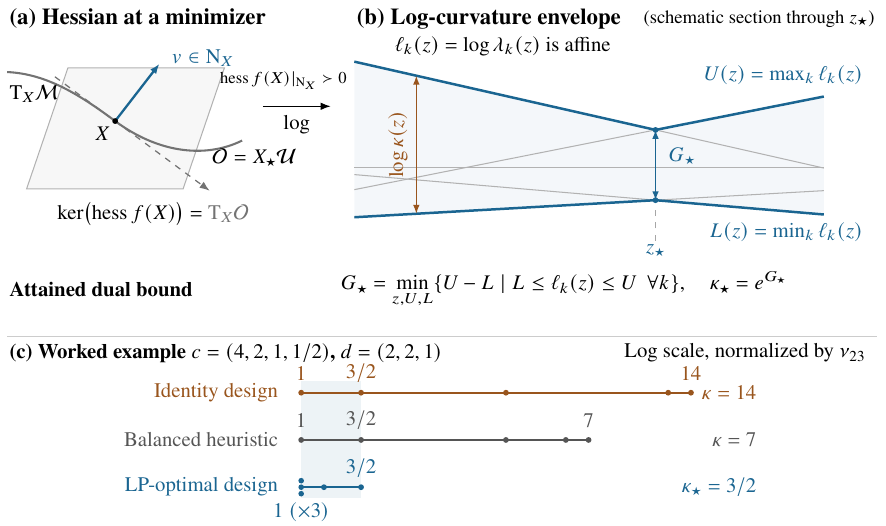}
    \caption{Minimizer geometry, logarithmic spread, and exact spectra. (a) At a minimizer, the orbit tangent is the Hessian kernel and the positive curvature acts on the metric-normal space. (b) Logarithms make the positive eigenvalues affine in the design, so that the LP minimizes the envelope spread at $z_\star$; the identity design's larger spread is shown for comparison. (c) The three designs' spectra, normalized by $\nu_{23}$, with the fixed interval $[1,3/2]$ shaded.}
    \label{fig:worked-spectral-design}
\end{figure}

In Figure~\ref{fig:worked-spectral-design}, under Theorem~\ref{thm:eigenvalue-Hf}, let $\mathcal M=\mathrm{Gq}(P,Q)$, $\mathcal O=X_\star\mathcal U$, and $H_X=\mathrm{hess}\,f(X)$ at $X\in\mathcal O$.
The metric-normal space $\mathrm N_X=(\mathrm T_X\mathcal O)^{\perp_{g_X^\beta}}\cap\mathrm T_X\mathcal M$ contains all positive Hessian modes, while $\ker H_X=\mathrm T_X\mathcal O$.
The positive eigenvalues $\lambda_k(z)$ of $H_X|_{\mathrm N_X}$ have affine logarithms $\ell_k(z)=\log\lambda_k(z)$. At an optimal equal-block design $z_\star$, their upper and lower envelopes have minimum spread $G_\star$.
This description includes isolated minimizers and zero optimal spread.
In the example, division by each design's $\nu_{23}$ places the positive spectra on a normalized logarithmic scale.
The shaded internal interval $[1,3/2]$ is fixed by the equal-block restriction.

    Fix $\tau=0$ and $x_1=x_2=0$ using Lemma~\ref{lem:LP-invariance}.
    The choice $x_3=-\log3$, $y_4=-\log4$ gives
    \begin{equation*}
        \varLambda_P=\mathrm{diag}(1,1,1/3),\qquad
        \varLambda_Q=\mathrm{diag}(1,1,1,1/4),\qquad \beta=1.
    \end{equation*}
    The leading three moduli of $Q$ do not enter the aligned spectrum and are set to $1$.
    In the row order of~\eqref{eq:worked-rows}, the positive eigenvalues are $(9/4,\ 3/2,\ 7/4,\ 3/2,\ 3/2)$.
    The feasible bounds $\underline\lambda=\log(3/2)$ and $\overline\lambda=\log(9/4)$ attain the lower bound.
    The \emph{identity design} in Remark~\ref{rem:special-cases} gives $(3/4,1/2,7,6,2)$ and condition number $14$.
    The \emph{balanced heuristic} uses $\varLambda_P=\mathrm{diag}(2,2,1)$ and $\varLambda_Q=\mathrm{diag}(4,2,1,1/2)$.
    Its positive eigenvalues are $(3/8,\ 1/4,\ 7/4,\ 3/2,\ 1)$, giving condition number $7$.
    The optimal design has
    \begin{equation*}
        G_\star=\log(3/2),\qquad \kappa_\star=3/2.
    \end{equation*}

    The dual vectors $u=(1,0,0,0,0)$ and $v=(0,1,0,0,0)$ are feasible in~\eqref{eq:dual-spread}.
    Both are nonnegative with unit sum, and $\xi_1=\xi_2$ gives the remaining dual equality.
    The primal and dual values both equal $\log(3/2)$.
    Rows $1$ and $2$ attain the upper and lower bounds, respectively, verifying complementary slackness.
    The sixth tangent mode is the orbit direction associated with $d_1=d_2$.
    The scalar design on that block makes the positive condition number constant along the orbit by Corollary~\ref{cor:diagonal-block}.

\FloatBarrier
\section{Implementation and Comparison with Preconditioned Metrics}\label{sec:metric-comparison}
The optimized constraint weights determine a realizable \textsf{PC} formulation and, through transport, a metric on the original manifold.
The implementation consists of four steps.
\begin{enumerate}
\item Form $\mathfrak C,\mathfrak D$ and ordered canonical decompositions. Drop zero modes and find equal-$d$ blocks.
\item Solve the spectral LP, fixing $\beta=1$ if desired. Recover $\alpha_i=e^{x_i}$ and $\gamma_j=e^{y_j}$, and choose positive leading moduli of $Q$.
\item Choose orthogonal bases, form $P,Q$, and construct $\breve P,\breve Q$ by~\eqref{eq:breve-construction}. Set $A,B$ by~\eqref{eq:equivalent-replacement}.
\item Pull an initial feasible $Y_0$ back to $X_0=\breve Q^{-1}Y_0\breve P$, apply a Riemannian solver with~\eqref{eq:metric-canonical}, and map the result back through $Y=\breve QX\breve P^{-1}$.
\end{enumerate}
The dense decompositions cost $O(n^3+p^3)$ and already identify minimizing eigenspaces for these quadratic problems.
Weights from an exact LP solution attain the optimal local condition number within the paired equal-block family.

The induced metric on the original manifold is obtained by transport.
For tangent vectors $U,V$ at $Y$, define
\begin{equation*}
    \widehat g_Y(U,V)=g^\beta_{\theta^{-1}(Y)}
    \bigl(\breve Q^{-1}U\breve P,\breve Q^{-1}V\breve P\bigr).
\end{equation*}
The map $\theta$ is an isometry, so that the Hessian spectra agree at corresponding critical points.
The LP minimizes the condition number over this family of transported metrics.

For comparison on $\mathrm{St}(p,n)$, we use the Lagrangian \textsf{PM} (\textsf{PM-Lag}) of \citet{mishra2016riemannian} and the left--right \textsf{PM} (\textsf{PM-LR}) of \citet{gao2025optimization}.
Table~\ref{tab:modal-comparison} compares their aligned modes with the \textsf{PC} design.
The table uses $\alpha_i=[\varLambda_P]_i$, $\gamma_j=[\varLambda_Q]_j$, $\beta=1$, and $s_i=\sqrt{(d_i/c_i)^2+\varepsilon}$ with $\varepsilon>0$.
The metric definitions, admissibility conditions, spectra, and proofs appear in the \prooflocation, Sections~\ref{sec:metric-Lagrangian} and~\ref{sec:metric-LR}.

\begin{table}[htbp]
    \centering
    \caption{Aligned Stiefel Hessian modes for the three design families.}
    \label{tab:modal-comparison}
    \renewcommand{\arraystretch}{1.8}
    \begin{tabular}{@{}lll@{}}
        \toprule
        Family & Internal mode, $i<j\le p$ & External mode, $i\le p<j$ \\
        \midrule
        \textsf{PC} & $\displaystyle\frac{(c_i-c_j)(d_i-d_j)}{\alpha_i\alpha_jc_ic_j}$ & $\displaystyle\frac{2d_i\gamma_j(c_i-c_j)}{\alpha_ic_ic_j}$ \\
        \textsf{PM-Lag} & $\displaystyle\frac{(c_i-c_j)(d_i-d_j)}{c_id_i+c_jd_j-\omega(c_id_j+c_jd_i)}$ & $\displaystyle\frac{c_i-c_j}{c_i-\omega c_j}$ \\
        \textsf{PM-LR} & $\displaystyle\frac{2(c_i-c_j)(d_i-d_j)}{c_is_i+c_js_j}$ & $\displaystyle\frac{2d_i(c_i-c_j)}{c_is_i}$ \\
        \bottomrule
    \end{tabular}
\end{table}

The \textsf{PC} weights enter the internal modes through $\alpha_i\alpha_j$ and the external modes through $\gamma_j/\alpha_i$.
Logarithms convert these coupled products and ratios into the affine LP rows.
The parameters $\omega$ and $\varepsilon$ act on the modal weights through the denominators shown in Table~\ref{tab:modal-comparison}.
The Lagrangian ambient inverse requires $0\le\omega<\omega_\star(Y)$, which at an aligned minimizer gives $\omega_\star=c_p/c_1$.
When the $c_i$ are distinct and $p>1$, the formal flattening limit $\omega\to1$ lies beyond that interval.
The LP optimum applies to the transported \textsf{PC} metric family. Comparisons with other metrics also depend on their modal weights and admissible parameter ranges.

Both \textsf{PM} constructions require an inverse metric action and a tangency solve.
The \textsf{PC} formulation uses the explicit tangent projector following~\eqref{eq:metric-canonical}.

\section{Numerical Experiments}
\label{sec:numerical}
We compare predicted positive spectra with returned-point diagnostics and track Riemannian conjugate gradient (R-CG) progress from a common initial point.
We test eight $(n,p)=(100,20)$ instances, each with one seed and starting point.
Total time includes preprocessing.

\subsection{Experimental Setup}
\label{sec:experimental-setup}
All runs use MATLAB R2025a and Manopt v8.0~\citep{boumal2014manopt} on an Apple M2 Max MacBook Pro with 12 cores and 32\,GB RAM.
The manifold implementation is available at \url{https://github.com/pxxyyz/GQmanifold}.
The computational archive records executed sources, solver settings, initializations, random states, and numerical diagnostics.

We use $C=\texttt{gallery('lehmer',n)}$.
The four choices of $D$ are $I_p$, $\mathrm{diag}(1,\ldots,p)$, $\mathrm{tridiag}(-1,2,-1)$, and \texttt{gallery('lehmer',p)}, labeled eigenvalue, Brockett, tridiagonal, and quadratic, respectively.
We compare the original \textsf{St}/\textsf{Sp} formulation with the balanced design \textsf{PC-bal} ($A=B=I$) and the LP design \textsf{PC-opt}.
For \textsf{St}, we also test \textsf{PM-Lag}, using~\eqref{eq:omega-adaptive} with $\eta=0.8$, and \textsf{PM-LR}, with $\varepsilon=10^{-3}$.
The original formulations use the Stiefel factory and the symplectic geometry of~\citet{gao2021riemannian,gao2024symplectic}, respectively.
The \textsf{PM} methods retain the Stiefel retraction and transport. Other methods use those supplied by their manifold factories.

We solve the design LP with MATLAB \texttt{linprog}.
With $\beta=1$, a second stage uses \texttt{quadprog} to minimize $\|x\|_2^2+\|y\|_2^2$ subject to spread at most $\widehat G+\epsilon\max\{1,|\widehat G|\}$, where $\widehat G$ is the numerical first-stage optimum.
The second stage accepts $\epsilon=10^{-9}$ on all eight instances.
The resulting \textsf{PC-opt} weights approximate the optimal design, and their spread may exceed the exact optimum $G_\star$.
The archive records accepted slacks, exit flags, and feasibility and optimality residuals.

Each driver uses random seed zero and reuses its initial point across the four choices of $D$.
Stiefel initialization is \texttt{orth(randn(n,p))}. Symplectic initialization follows the block-prefix construction of~\citet[Sec.~6.1]{gao2021riemannian}, with the executed formula recorded in the archive.
Each \textsf{PC} formulation starts from the inverse image of the same initial point under $\theta$.
All runs use Manopt's Hestenes--Stiefel update and backtracking line search, a gradient tolerance of $10^{-6}$, and an iteration cap of $10^4$, retaining the other stopping defaults.
Since the gradient norms depend on the metric, the common tolerance permits different stationarity accuracies in the original coordinates.

Table~\ref{tab:quadratic} reports objective value (\texttt{fval}), Riemannian gradient norm (\texttt{gradf}), feasibility error $\|X^{(P,Q)}X-I_p\|_{\mathrm F}$ (\texttt{feasi}), iterations, preprocessing (\texttt{setup}), and total time in seconds.
Preprocessing includes initial-point generation, required spectral work, metric or constraint construction, both design stages, and options setup.
We charge this work once to each standalone method and add the outer solver wall time to obtain total time.
Given-data construction, reference optima, condition diagnostics, plotting, and result file output are excluded.
Plot endpoints equal total time. Intermediate points use shifted internal solver clocks.

The predicted aligned condition number is \texttt{condhat}.
The diagnostic \texttt{cond} uses \texttt{hessian\-spectrum}, retaining real parts above $10^{-12}\max\{1,\max_i|\lambda_i|\}$.
For \textsf{PC}, the supplied operator is the Riemannian Hessian.
For \textsf{PM-Lag} and \textsf{PM-LR}, it is the preconditioned Lagrangian operator~\eqref{eq:H-hess}, which agrees with the Hessian at critical points but omits metric-derivative terms elsewhere.
For the original \textsf{Sp} cost, the aligned formula does not apply, and \texttt{condhat} is \texttt{NaN}.

\subsection{Results and Interpretation}
\begin{figure}[p]
    \begin{minipage}[t][0.98\textheight][s]{\linewidth}
        \centering
        \includegraphics[width=0.680\linewidth]{\datapath 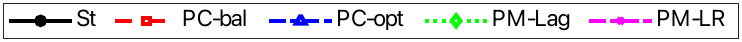}
        \par\vfill
        \begin{minipage}{\figwidth}
            \centering
            \includegraphics[width=\linewidth]{\datapath 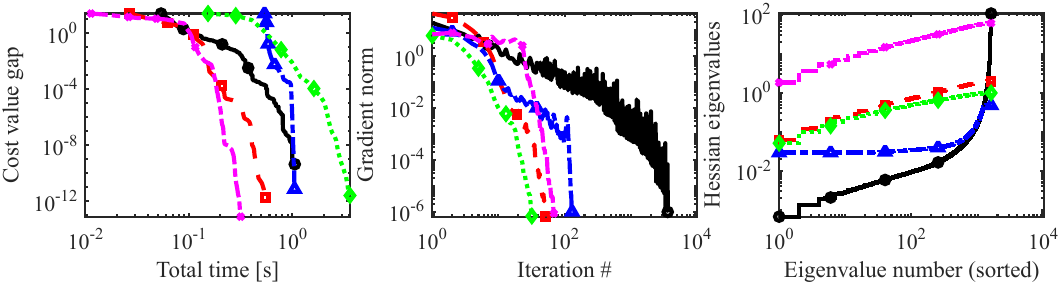}
            \par\footnotesize(a) \textbf{Symmetric eigenvalue problem.}

        \end{minipage}

        \par\vfill

        \begin{minipage}{\figwidth}
            \centering
            \includegraphics[width=\linewidth]{\datapath 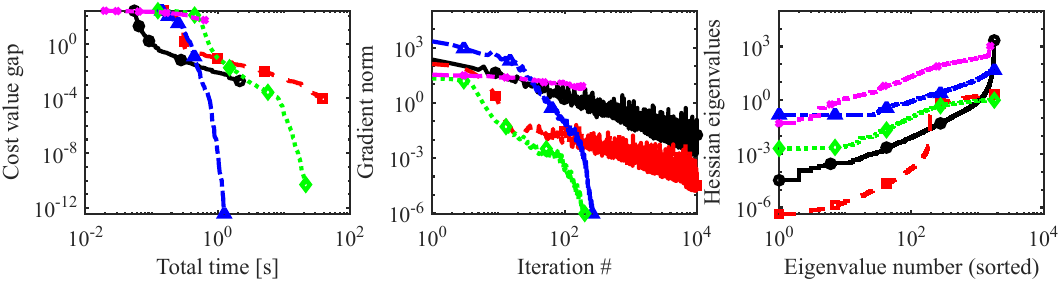}
            \par\footnotesize(b) \textbf{Symmetric Brockett problem.}

        \end{minipage}

        \par\vfill

        \begin{minipage}{\figwidth}
            \centering
            \includegraphics[width=\linewidth]{\datapath 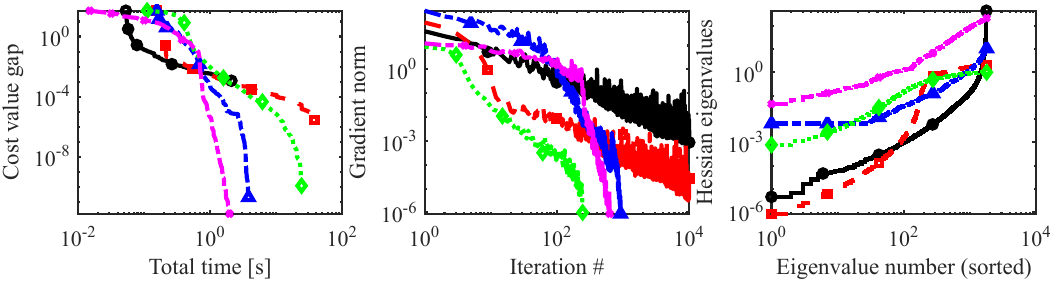}
            \par\footnotesize(c) \textbf{Symmetric tridiagonal problem.}

        \end{minipage}

        \par\vfill

        \begin{minipage}{\figwidth}
            \centering
            \includegraphics[width=\linewidth]{\datapath 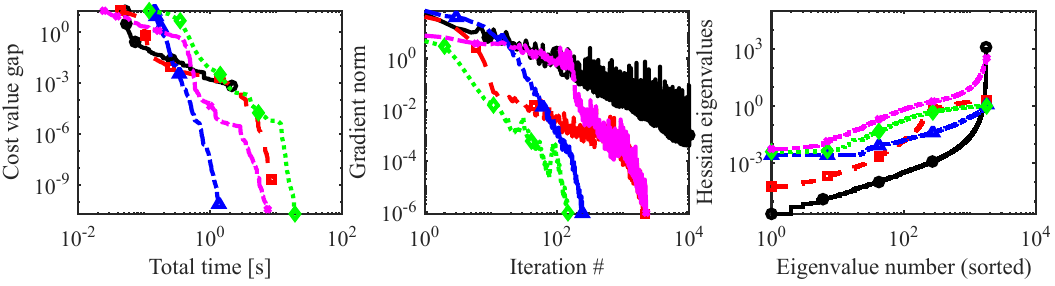}
            \par\footnotesize(d) \textbf{Symmetric quadratic problem.}

        \end{minipage}

        \par\vfill
        \caption{Convergence histories and returned-point diagnostic spectra for the four symmetric instances.}
        \label{fig:convergence_quadratic_symm}
    \end{minipage}
    \end{figure}

\begin{figure}[p]
    \begin{minipage}[t][0.98\textheight][s]{\linewidth}
        \centering
        \includegraphics[width=0.389\linewidth]{\datapath 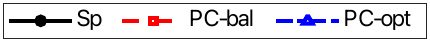}
        \par\vfill
        \begin{minipage}{\figwidth}
            \centering
            \includegraphics[width=\linewidth]{\datapath 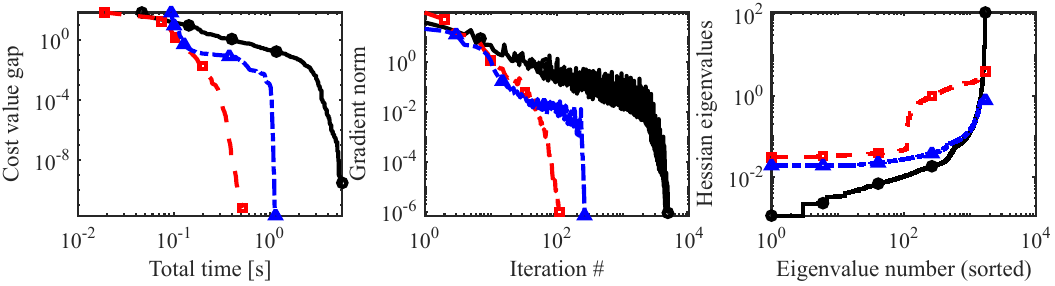}
            \par\footnotesize(a) \textbf{Skew-symmetric eigenvalue problem.}

        \end{minipage}

        \par\vfill

        \begin{minipage}{\figwidth}
            \centering
            \includegraphics[width=\linewidth]{\datapath 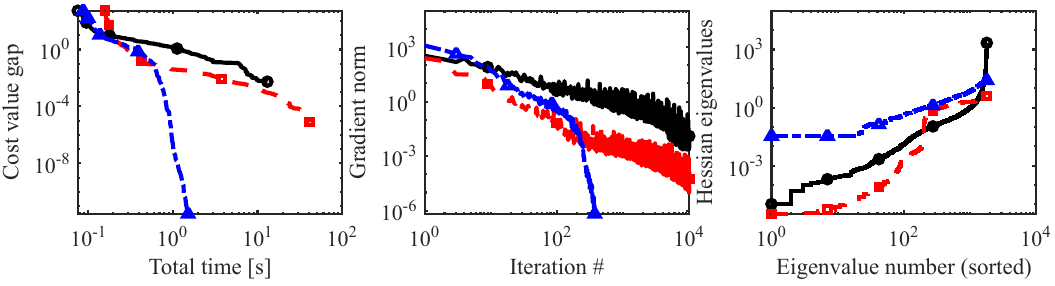}
            \par\footnotesize(b) \textbf{Skew-symmetric Brockett problem.}

        \end{minipage}

        \par\vfill

        \begin{minipage}{\figwidth}
            \centering
            \includegraphics[width=\linewidth]{\datapath 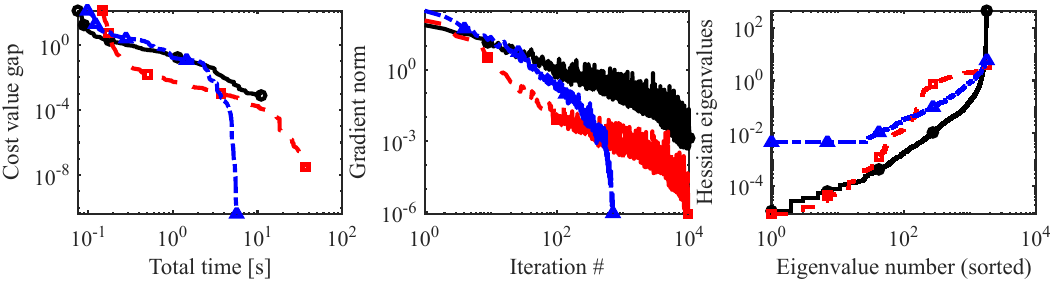}
            \par\footnotesize(c) \textbf{Skew-symmetric tridiagonal problem.}

        \end{minipage}

        \par\vfill

        \begin{minipage}{\figwidth}
            \centering
            \includegraphics[width=\linewidth]{\datapath 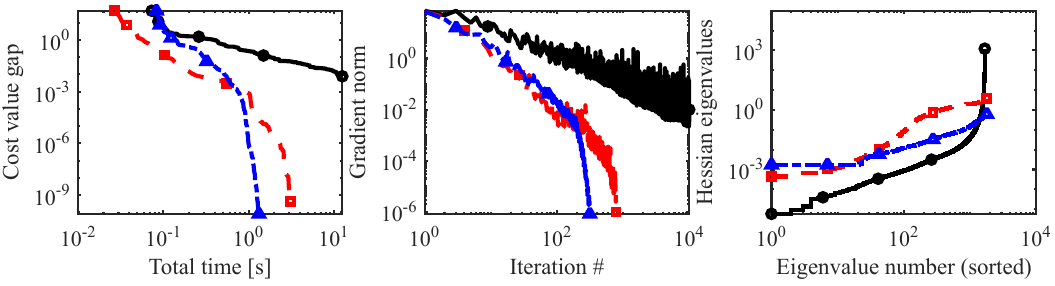}
            \par\footnotesize(d) \textbf{Skew-symmetric quadratic problem.}

        \end{minipage}

        \par\vfill
        \caption{Convergence histories and returned-point diagnostic spectra for the four skew-symmetric instances.}
        \label{fig:convergence_quadratic_skew}
    \end{minipage}
    \end{figure}

Figures~\ref{fig:convergence_quadratic_symm} and~\ref{fig:convergence_quadratic_skew} display, from top to bottom, the eigenvalue, Brockett, tridiagonal, and quadratic instances.
The three columns show the cost gap against total time, gradient norm against iterations, and sorted retained positive real parts of diagnostic eigenvalues at the returned point.
All axes are logarithmic, and the optimal value is given by Proposition~\ref{prop:critical-point}.
Table~\ref{tab:quadratic} reports the aligned prediction \texttt{condhat} alongside the returned-point diagnostic \texttt{cond}.

Some runs terminating at the iteration cap or minimum step size have diagnostic eigenvalues with negative real parts.
At these nonstationary outputs, \texttt{cond} uses only the retained positive spectrum and may differ from \texttt{condhat}.
Negative real parts at a nonstationary point do not determine the type of a nearby critical orbit.

On the tested orthogonality-constrained problems, the LP design reduces the predicted local condition number by several orders of magnitude relative to the original formulation.
The numerical LP design reaches the gradient tolerance on all eight tested instances. Its returned positive-spectrum condition numbers closely match the aligned predictions.
On the symplectic instances, these condition numbers are substantially smaller than those of \textsf{Sp}, which reaches the iteration cap on the Brockett, tridiagonal, and quadratic problems.
On the eigenvalue problems, \textsf{PC-opt} has a smaller positive condition number than \textsf{PC-bal} but requires more iterations and longer total time.
Improved local conditioning alone does not guarantee faster convergence for the initializations considered.

On the symmetric instances, \textsf{PM-Lag} reaches the gradient tolerance in fewer iterations than \textsf{PC-opt}, but its total time is longer.
The comparison distinguishes iteration counts from total computational cost.
For \textsf{PM-LR}, the Brockett run terminates with a zero step while its gradient norm remains above tolerance. The remaining runs reach the gradient tolerance.

These fixed-size experiments support the spectral predictions and illustrate solver progress under metric-dependent stopping norms.
Assessing robustness and scalability requires repeated instances across dimensions and a common stationarity residual in the original variables.

\clearpage
\begin{table}[ht]
    \centering
    \footnotesize
    \caption{Validation on eight quadratic instances.}
    \label{tab:quadratic}
    \pgfplotstableset{
        col sep=comma,
        columns={method,fval,gradf,feasi,iter,setup,time,cond,cond hat},
        every head row/.style={before row=\toprule, after row=\midrule},
        every last row/.style={after row=\bottomrule},
        columns/Row/.style={column name={}, string type, column type={>{\sffamily}w{r}{1cm}}},
        columns/method/.style={column name={}, string type, column type={>{\sffamily}w{r}{1.05cm}}},
        columns/fval/.style={column name=\texttt{fval}, sci, precision=2, column type=w{r}{1.4cm}},
        columns/gradf/.style={column name=\texttt{gradf}, sci, precision=2, column type=w{r}{1.4cm}},
        columns/feasi/.style={column name=\texttt{feasi}, sci, precision=2, column type=w{r}{1.4cm}},
        columns/iter/.style={column name=\texttt{iter}, fixed, precision=0, column type=w{r}{1cm}},
        columns/time/.style={column name=\texttt{total}, fixed, precision=2, column type=w{r}{1cm}},
        columns/setup/.style={column name=\texttt{setup}, sci, precision=2, column type=w{r}{1.45cm}},
        columns/cond/.style={column name=\texttt{cond}, sci, precision=2, column type=w{r}{1.2cm}},
        columns/cond hat/.style={column name=\texttt{condhat}, sci, precision=2, column type=w{r}{1.2cm}},
        column type=r,
        every column/.append style={font=\footnotesize},
    }

    \pgfplotstabletypeset[
        every first row/.style={
            before row={
                \multicolumn{6}{l}{
                \textbf{Symmetric eigenvalue problem}
                } \\
            },
            after row=,
        },
    ]{\datapath tab-symm-eigenvalue.csv}

    \pgfplotstabletypeset[
        every first row/.style={
            before row={
                \multicolumn{6}{l}{
                \textbf{Symmetric Brockett problem}
                } \\
            },
            after row=,
        },
        every head row/.style={output empty row},
    ]{\datapath tab-symm-brockett.csv}

    \pgfplotstabletypeset[
        every first row/.style={
            before row={
                \multicolumn{6}{l}{
                \textbf{Symmetric tridiagonal problem}
                } \\
            },
            after row=,
        },
        every head row/.style={output empty row},
    ]{\datapath tab-symm-tridiagonal.csv}

    \pgfplotstabletypeset[
        every first row/.style={
            before row={
                \multicolumn{6}{l}{
                \textbf{Symmetric quadratic problem}
                } \\
            },
            after row=,
        },
        every head row/.style={output empty row},
    ]{\datapath tab-symm-quadratic.csv}

    \pgfplotstabletypeset[
        every first row/.style={
            before row={
                \multicolumn{6}{l}{
                \textbf{Skew-symmetric eigenvalue problem}
                } \\
            },
            after row=,
        },
        every head row/.style={output empty row},
    ]{\datapath tab-skew-eigenvalue.csv}

    \pgfplotstabletypeset[
        every first row/.style={
            before row={
                \multicolumn{6}{l}{
                \textbf{Skew-symmetric Brockett problem}
                } \\
            },
            after row=,
        },
        every head row/.style={output empty row},
    ]{\datapath tab-skew-brockett.csv}

    \pgfplotstabletypeset[
        every first row/.style={
            before row={
                \multicolumn{6}{l}{
                \textbf{Skew-symmetric tridiagonal problem}
                } \\
            },
            after row=,
        },
        every head row/.style={output empty row},
    ]{\datapath tab-skew-tridiagonal.csv}

    \pgfplotstabletypeset[
        every first row/.style={
            before row={
                \multicolumn{6}{l}{
                \textbf{Skew-symmetric quadratic problem}
                } \\
            },
            after row=,
        },
        every head row/.style={output empty row},
    ]{\datapath tab-skew-quadratic.csv}
    \end{table}

\section{Conclusions and Perspectives}
\label{sec:conclusion}
Reshaping the constraints provides design variables for preconditioning quadratic Stiefel and symplectic problems.
The paired congruence construction and equal-block restriction specify the admissible constraint family.
At the quotient minimizer, simultaneous diagonalization reduces positive Hessian condition-number minimization over this family exactly to an LP.
The attained logarithmic spread $G_\star$ gives the optimal condition number $e^{G_\star}$ within this family, and LP strong duality gives a matching structural lower bound.

The exact local optimum provides a reference for designs based on approximate or reusable spectral information.
The eight fixed-size experiments support the spectral formulas and examine R-CG progress, with preprocessing included in total time.
Extension to nonlinear objectives or congruences beyond the paired construction requires new curvature formulas.

\FloatBarrier
\nocite{hardy1952inequalities}
\paragraph{Funding.} This work was supported by the National Key R\&D Program of China (grant No. 2023YFA1009300). 
    BG was supported by the National Natural Science Foundation of China (grant No. 12288201). 
    ZP was supported by the China Postdoctoral Science Foundation (grant No. 2024M763472).
\bibliographystyle{plainnat}
\bibliography{references}

@misc{peng2026generalized,
  title = {Generalized quadratic matrix manifolds: {R}iemannian geometry and computations},
  author = {Peng, Zhen and Gao, Bin},
  howpublished = {\emph{arXiv preprint arXiv:2610.28159}},
  year = {2026}
}

@book{absil2008optimization,
  title = {Optimization algorithms on matrix manifolds},
  author = {Absil, P.-A. and Mahony, Robert and Sepulchre, Rodolphe},
  year = {2008},
  address = {Princeton, NJ},
  publisher = {Princeton University Press}
}

@inproceedings{absil2013extrinsic,
  title = {An extrinsic look at the {R}iemannian {H}essian},
  author = {Absil, P.-A. and Mahony, Robert and Trumpf, Jochen},
  booktitle = {Geometric Science of Information},
  series = {Lecture Notes in Computer Science},
  volume = {8085},
  pages = {361--368},
  year = {2013},
  publisher = {Springer},
  address = {Berlin, Heidelberg}
}

@article{boumal2014manopt,
  title = {Manopt, a {MATLAB} toolbox for optimization on manifolds},
  author = {Boumal, Nicolas and Mishra, Bamdev and Absil, P.-A. and Sepulchre, Rodolphe},
  journal = {J. Mach. Learn. Res.},
  shortjournal = {J. Mach. Learn. Res.},
  volume = {15},
  number = {42},
  pages = {1455--1459},
  year = {2014},
  publisher = {JMLR.org}
}

@book{boumal2023introduction,
  title     = {An introduction to optimization on smooth manifolds},
  author    = {Boumal, Nicolas},
  year      = {2023},
  address   = {Cambridge, UK},
  publisher = {Cambridge University Press}
}

@article{gao2021riemannian,
  title = {{R}iemannian optimization on the symplectic {S}tiefel manifold},
  author = {Gao, Bin and Son, Nguyen Thanh and Absil, P.-A. and Stykel, Tatjana},
  journal = {SIAM J. Optim.},
  shortjournal = {SIAM J. Optim.},
  volume = {31},
  number = {2},
  pages = {1546--1575},
  year = {2021},
  publisher = {SIAM}
}

@misc{gao2024symplectic,
  title   = {Symplectic {S}tiefel manifold: tractable metrics, second-order geometry and {N}ewton's methods},
  author  = {Gao, Bin and Son, Nguyen Thanh and Stykel, Tatjana},
  howpublished = {\emph{arXiv preprint arXiv:2406.14299}},
  year    = {2024}
}

@article{gao2025optimization,
  title     = {Optimization on product manifolds under a preconditioned metric},
  author    = {Gao, Bin and Peng, Renfeng and Yuan, {Ya-xiang}},
  journal   = {SIAM J. Matrix Anal. Appl.},
  shortjournal = {SIAM J. Matrix Anal. Appl.},
  volume    = {46},
  number    = {3},
  pages     = {1816--1845},
  year      = {2025},
  publisher = {SIAM}
}

@book{hardy1952inequalities,
  title={Inequalities},
  author={Hardy, Godfrey Harold and Littlewood, John Edensor and P{\'o}lya, George},
  year={1952},
  publisher={Cambridge University Press}
}

@article{liu2019quadratic,
  title     = {Quadratic optimization with orthogonality constraint: explicit {{\L}}ojasiewicz exponent and linear convergence of retraction-based line-search and stochastic variance-reduced gradient methods},
  author    = {Liu, Huikang and So, Anthony Man-Cho and Wu, Weijie},
  journal   = {Math. Program.},
  shortjournal = {Math. Program.},
  volume    = {178},
  number    = {1-2},
  pages     = {215--262},
  year      = {2019},
  publisher = {Springer}
}

@article{mishra2016riemannian,
  title = {{R}iemannian preconditioning},
  author = {Mishra, Bamdev and Sepulchre, Rodolphe},
  journal = {SIAM J. Optim.},
  shortjournal = {SIAM J. Optim.},
  volume = {26},
  number = {1},
  pages = {635--660},
  year = {2016},
  publisher = {SIAM}
}

@article{shustin2023riemannian,
  title = {{R}iemannian optimization with a preconditioning scheme on the generalized {S}tiefel manifold},
  author = {Shustin, Boris and Avron, Haim},
  journal = {J. Comput. Appl. Math.},
  shortjournal = {J. Comput. Appl. Math.},
  volume = {423},
  pages = {114953},
  year = {2023},
  publisher = {Elsevier}
}

@article{shustin2024faster,
  title   = {Faster randomized methods for orthogonality constrained problems},
  author  = {Shustin, Boris and Avron, Haim},
  journal = {J. Mach. Learn. Res.},
  shortjournal = {J. Mach. Learn. Res.},
  volume  = {25},
  number  = {257},
  pages   = {1--59},
  year    = {2024}
}

@article{son2021computing,
  title = {Computing symplectic eigenpairs of symmetric positive-definite matrices via trace minimization and {R}iemannian optimization},
  author = {Son, Nguyen Thanh and Absil, P.-A. and Gao, Bin and Stykel, Tatjana},
  journal = {SIAM J. Matrix Anal. Appl.},
  shortjournal = {SIAM J. Matrix Anal. Appl.},
  volume = {42},
  number = {4},
  pages = {1732--1757},
  year = {2021},
  publisher = {SIAM}
}

@misc{son2025brockett,
  title   = {Brockett cost function for symplectic eigenvalues},
  author  = {Son, Nguyen Thanh},
  howpublished = {\emph{arXiv preprint arXiv:2506.07560}},
  year    = {2025}
}

@article{edelman1998geometry,
  title={The geometry of algorithms with orthogonality constraints},
  author={Edelman, Alan and Arias, Tom{\'a}s A and Smith, Steven T},
  journal={SIAM J. Matrix Anal. Appl.},
  shortjournal = {SIAM J. Matrix Anal. Appl.},
  volume={20},
  number={2},
  pages={303--353},
  year={1998},
  publisher={SIAM}
}

@misc{dougan2025optimal,
  title = {Optimal preconditioning is a geodesically convex optimization problem},
  author = {Do{\u{g}}an, M. Levent and Erg{\"u}r, Alperen and Tsigaridas, Elias},
  howpublished = {\emph{arXiv preprint arXiv:2512.06618}},
  year = {2025}
}

@article{braatz1994minimizing,
  title={Minimizing the {Euclidean} condition number},
  author={Braatz, Richard D and Morari, Manfred},
  journal={SIAM J. Control Optim.},
  shortjournal = {SIAM J. Control Optim.},
  volume={32},
  number={6},
  pages={1763--1768},
  year={1994},
  publisher={SIAM}
}

@article{qu2025optimal,
  title={Optimal diagonal preconditioning},
  author={Qu, Zhaonan and Gao, Wenzhi and Hinder, Oliver and Ye, Yinyu and Zhou, Zhengyuan},
  journal={Oper. Res.},
  shortjournal = {Oper. Res.},
  volume={73},
  number={3},
  pages={1479--1495},
  year={2025},
  publisher={INFORMS}
}

@inproceedings{jambulapati2023structured,
  title = {Structured semidefinite programming for recovering structured preconditioners},
  author = {Jambulapati, Arun and Li, Jerry and Musco, Christopher and Shiragur, Kirankumar and Sidford, Aaron and Tian, Kevin},
  booktitle = {Advances in Neural Information Processing Systems},
  volume = {36},
  pages = {16991--17007},
  year = {2023},
  publisher = {Curran Associates, Inc.}
}

@book{boyd1994linear,
  title={Linear matrix inequalities in system and control theory},
  author={Boyd, Stephen and El Ghaoui, Laurent and Feron, Eric and Balakrishnan, Venkataramanan},
  year={1994},
  publisher={SIAM}
}

@article{boyd2007tutorial,
  title={A tutorial on geometric programming},
  author={Boyd, Stephen and Kim, Seung-Jean and Vandenberghe, Lieven and Hassibi, Arash},
  journal={Optim. Eng.},
  shortjournal = {Optim. Eng.},
  volume={8},
  number={1},
  pages={67--127},
  year={2007},
  publisher={Springer}
}

@book{vandenberghe2004convex,
  title = {Convex optimization},
  author = {Boyd, Stephen and Vandenberghe, Lieven},
  year = {2004},
  publisher = {Cambridge University Press},
  address = {Cambridge, UK}
}
\clearpage
\appendix
\renewcommand{\theHsection}{appendix.\Alph{section}}
\normalsize

\numberwithin{equation}{section}
\numberwithin{lemma}{section}
\numberwithin{proposition}{section}
\numberwithin{remark}{section}
\section{Proofs of the Main Results}\label{app:proofs}

\subsection{Proof of Lemma~\ref{lem:critical-point}}
\begin{namedproof}{Proof.}
At a critical point, the normal-space condition gives $AXB=QXP^{-1}X^\top AXB$ and~\eqref{eq:critical-point-eq}.
        Invertibility of $B$ yields $AX=QXP^{-1}X^\top AX$, and hence $X_\perp^\top AX=0$.
        The selected and complementary blocks therefore satisfy
        \begin{equation*}
            E^{-1}Q^{-1}AE
            =\mathrm{blkdiag}(P^{-1}X^\top AX,\,O^{-1}X_\perp^\top AX_\perp),
        \end{equation*}
        which proves the complementary-spectrum claim.
        The eigenvalues of the symmetric congruence $O^{-\top}X_\perp^\top AX_\perp O^{-1}$ generally differ from those of the complementary pencil.
\end{namedproof}

\subsection{Proof of Proposition~\ref{prop:Rayleigh-quotient}}
\begin{namedproof}{Proof.}
Lemma~\ref{lem:critical-point} gives $X_\perp^\top AX=X_\perp^\top QX=0$.
        Substituting~\eqref{eq:normalized-tangent} into~\eqref{eq:Rayleigh-quotient}, with $u=\mathrm{vec}(\widehat\Omega)$ and $v=\mathrm{vec}(V)$, gives
        \begin{equation*}
            q(Z)=\frac{2\beta^{-1}u^\top H_{\mathrm{int}}u+2v^\top H_{\mathrm{ext}}v}
            {u^\top u+v^\top v},\qquad u=\Pi_\varrho u.
        \end{equation*}
        The restriction to $\mathrm{range}(\Pi_\varrho)$ gives the intrinsic internal operator $2\Pi_\varrho H_{\mathrm{int}}\Pi_\varrho/\beta$.
        Relative to this subspace and its orthogonal complement, $2\Pi_\varrho H_{\mathrm{int}}/\beta$ is block upper triangular.
        Its diagonal blocks are the intrinsic operator and a zero block. This proves the nonzero-spectrum identity in~\eqref{eq:Hf}.
\end{namedproof}

\subsection{Proof of Proposition~\ref{prop:commutation}}
\begin{namedproof}{Proof.}
Use the real canonical decompositions $Q = \varPhi_Q\varSigma_Q\varPhi_Q^\top$ and $P = \varPhi_P\varSigma_P\varPhi_P^\top$, with $\varSigma_Q,\varSigma_P$ as in~\eqref{eq:spectral-sigma}.
        Write the fixed data matrices as $\mathfrak C=\varPhi_{\mathfrak C}\varSigma_{\mathfrak C}\varPhi_{\mathfrak C}^{\top}$ and $\mathfrak D=\varPhi_{\mathfrak D}\varSigma_{\mathfrak D}\varPhi_{\mathfrak D}^{\top}$.
        Define
        \begin{equation}
            \breve{Q} := C^{-1/2}\varPhi_{\mathfrak{C}}\,(\varSigma_Q\varSigma_{\mathfrak{C}}^{-1})^{1/2}\,\varPhi_{Q}^{\top},
            \qquad
            \breve{P} := D^{1/2}\varPhi_{\mathfrak{D}}\,(\varSigma_P\varSigma_{\mathfrak{D}}^{-1})^{1/2}\,\varPhi_{P}^{\top}.
        \end{equation}
        Positive definiteness of the products under the square roots makes both factors well defined.
        Substitution gives~\eqref{eq:congruence} and
        \begin{equation*}
            A = \varPhi_{Q}\,\varSigma_Q\varSigma_{\mathfrak{C}}^{-1}\,\varPhi_{Q}^{\top},
            \qquad
            B = \varPhi_{P}\,\varSigma_P^{-1}\varSigma_{\mathfrak{D}}\,\varPhi_{P}^{\top},
        \end{equation*}
        which proves (a) and (c).
        Claim (b) follows from $A^{-1}Q = \varPhi_Q\varSigma_{\mathfrak{C}}\varPhi_Q^{\top}$ and $BP = \varPhi_P\varSigma_{\mathfrak{D}}\varPhi_P^{\top}$.
\end{namedproof}

\subsection{Proof of Proposition~\ref{prop:critical-point}}
\begin{namedproof}{Proof.}
We first classify critical points in the paired canonical frames of Proposition~\ref{prop:commutation} and then determine the minimizing assignment.
All roots in this proof are positive canonical modulus roots, in both the symmetric and skew cases.
By Lemma~\ref{lem:critical-point}, the range $\mathrm{range}(X)$ at a critical point is invariant under $Q^{-1}A$.
The simple values in $c$ force this range to be spanned by selected eigenlines for $\varrho=+1$ or selected invariant two-planes for $\varrho=-1$.
For an ordered tuple $\mathcal I$ of these modes, write
$X=[\varPhi_Q]_{:,\mathcal I}\varPi\varPhi_P^\top$ and put
\begin{equation*}
 F=[\varSigma_{\mathfrak C}]_{\mathcal I}^{-1}K_\varrho(P),
 \qquad \widehat{Z}=[\varSigma_Q]_{\mathcal I}^{1/2}\varPi\varSigma_P^{-1/2},
 \qquad \widehat{H}=\widehat{Z}^\top F\widehat{Z}.
\end{equation*}
The matrices $F$ and $\widehat{H}$ are SPD.
The constraint becomes $\widehat{Z}^\top K_\varrho(P)\widehat{Z}=K_\varrho(P)$, while stationarity gives
$\widehat{H}\varSigma_{\mathfrak D}=\varSigma_{\mathfrak D}\widehat{H}$.
Let $\mathcal C$ be the orthogonal centralizer of $\varSigma_{\mathfrak D}$ preserving $K_\varrho(P)$.
Its equal-$d$ blocks are $\mathrm O(m)$ in the symmetric case and $\mathrm U(m)$ in the skew case.

For $\varrho=+1$, $\widehat{Z}$ is orthogonal.
Diagonalizing $\widehat{H}$ within its equal-$d$ blocks and ordering the selected indices gives $\widehat{H}=W^\top FW$ for some $W\in\mathcal C$.
The matrix $\widehat{Z}W^\top$ is orthogonal and preserves the simple diagonal $F$. Its diagonal signs can be absorbed into $W$.
For $\varrho=-1$, identify each $J_2$ plane with a complex coordinate.
Within each equal-$d$ block, the real symmetric matrix $\widehat{H}$ commutes with the complex structure and therefore represents a Hermitian matrix.
Its positive eigenvalues equal its symplectic frequencies.
With $K=K_\varrho(P)$, symplecticity gives
\begin{equation*}
K^{-1}\widehat{H}=\widehat{Z}^{-1}(K^{-1}F)\widehat{Z}.
\end{equation*}
Unitary diagonalization therefore recovers the selected frequencies $1/c_i$.
Ordering the selected planes again gives $\widehat{H}=W^\top FW$ with $W\in\mathcal C$.
Set $\widetilde{S}=\widehat{Z}W^\top$.
The identities $\widetilde{S}^\top K\widetilde{S}=K$ and $\widetilde{S}^\top F\widetilde{S}=F$ imply $\widetilde{S}^{-1}(K^{-1}F)\widetilde{S}=K^{-1}F$.
The distinct frequencies force $\widetilde{S}$ to preserve each selected plane.
On each plane it commutes with $J_2$ and has the form $aI_2+bJ_2$.
Preservation of $F$ gives $a^2+b^2=1$, so that the block is a rotation and can be absorbed into $W$.
Every critical point therefore has the representation
\begin{equation*}
 X=[\varPhi_Q]_{:,\mathcal I}[\varSigma_Q]_{\mathcal I}^{-1/2}
 W\varSigma_P^{1/2}\varPhi_P^\top,\qquad W\in\mathcal C.
\end{equation*}
Conversely, direct substitution verifies stationarity and feasibility.
The corresponding objective value is $\mathrm{tr}([\varSigma_{\mathfrak C}]_{\mathcal I}^{-1}\varSigma_{\mathfrak D})$.
As $A,B$ are positive definite, the objective is coercive and attains its minimum on the nonempty closed feasible set.
The rearrangement inequality~\citep[Thm.~368]{hardy1952inequalities} assigns the largest $c_i$ in the order of decreasing $d_i$.
The simple values in $c$ make this assignment unique modulo equal-$d$ permutations. The stated minimum and minimizing orbit follow.
\end{namedproof}

\subsection{Proof of Proposition~\ref{prop:critical-manifold}}
\begin{namedproof}{Proof.}
\emph{Group action.}
For $U$ in~\eqref{eq:gauge-group}, set
$R=\varSigma_P^{1/2}\varPhi_P^\top U\varPhi_P\varSigma_P^{-1/2}$.
The defining identities of $U$ give
$R^\top K_pR=K_p$ and $RD_0R^\top=D_0$, where
$K_p=K_\varrho(P)$ and $D_0=\mathrm{diag}(d_i)\otimes I_\nu$.
The matrix $S=D_0^{-1/2}RD_0^{1/2}$ is orthogonal and centralizes $D_0^{1/2}K_pD_0^{1/2}$.
The distinct canonical moduli are the distinct values in $d$. Thus $S$ preserves their blocks and commutes with $D_0$.
It follows that $R=S$ is orthogonal.
The correspondence gives an isomorphism from $\mathcal U$ onto the centralizer $\mathcal C_D$ in Section~\ref{sec:spectrum}.
Thus $\mathcal U$ is compact and isomorphic to the stated product of orthogonal or unitary groups.
The preceding critical-point classification gives the compact embedded minimizing orbit $X_\star\mathcal U$ with dimension~\eqref{eq:dim-centralizer}.
Full column rank makes the action free, as $XU=X$ implies $U=I$. Compactness makes it proper.
For $X=X_\star U_0$, the identity $U_0\mathcal U=\mathcal U$ gives $X\mathcal U=X_\star\mathcal U$.

For every positive design, the canonical metric has the equivalent form
\begin{equation}\label{eq:metric-projector-form}
 M_X^\beta=\beta Q^\top XX^\top Q+I-X(X^\top X)^{-1}X^\top.
\end{equation}
Indeed, the block rows of $E^{-1}$ are $P^{-1}X^\top Q$ and $O^{-1}X_\perp^\top Q$.
Substitution in~\eqref{eq:metric-canonical} gives the first term and $Q^\top X_\perp X_\perp^\top Q$.
The columns of $Q^\top X_\perp$ form an orthonormal basis for $\ker X^\top$. The sum of their outer products gives the displayed projector.

\emph{Direct mode calculation.}
Write $A_\alpha=\mathrm{diag}(\alpha_i)\otimes I_\nu$ and
$A_\gamma=\mathrm{diag}(\gamma_j)\otimes I_\nu$.
Thus $\varSigma_P^{1/2}=A_\alpha^{1/2}$ and
$\varSigma_Q^{1/2}=A_\gamma^{1/2}$.
Substitution of the paired formulas for $A,B$ gives the fixed constraint and cost stated after~\eqref{eq:common-reference}.
At $W_\star=[I_p,0]^\top$, write a reference tangent as
$\Delta=[\Omega^\top,V^\top]^\top$, with
$\Omega^\top K_p+K_p\Omega=0$.
The corresponding tangent vector in the original realization is
$Z=\varPhi_QA_\gamma^{-1/2}\Delta A_\alpha^{1/2}\varPhi_P^\top$.
Formula~\eqref{eq:metric-projector-form} gives
\begin{equation}\label{eq:aligned-metric-pair}
 g_{X_\star}^\beta(Z,Z)
 =\beta\|A_\alpha^{1/2}\Omega A_\alpha^{1/2}\|_{\mathrm F}^2
 +\mathrm{tr}(V^\top[A_\gamma]_{>p}^{-1}VA_\alpha).
\end{equation}
The second differential of the fixed reference cost is
\begin{equation}\label{eq:reference-second-form}
 2\mathrm{tr}(\Delta^\top C_0\Delta D_0-\Delta^\top K_n\Delta M),
 \qquad M=K_p^{-1}W_\star^\top C_0W_\star D_0.
\end{equation}

For $\varrho=+1$, use Frobenius-unit internal matrices
$(e_ie_j^\top-e_je_i^\top)/\sqrt2$ and unit external entries.
Their metric coefficients in~\eqref{eq:aligned-metric-pair} are
$\beta\alpha_i\alpha_j$ and $\alpha_i/\gamma_j$.
Their curvature coefficients in~\eqref{eq:reference-second-form} are
$(d_i-d_j)(c_i-c_j)/(c_ic_j)$ and $2d_i(c_i-c_j)/(c_ic_j)$.
Distinct coordinates are orthogonal for both forms.

For $\varrho=-1$, an off-diagonal internal block $\Omega_{ij}$ determines
$\Omega_{ji}=J_2\Omega_{ij}^\top J_2$.
The matrices $I_2,J_2$ commute with $J_2$ and give the two minus modes.
The matrices $\mathrm{diag}(1,-1)$ and
$\left[\begin{smallmatrix}0&1\\1&0\end{smallmatrix}\right]$ anticommute with $J_2$ and give the two plus modes.
Divide each resulting two-block matrix by $2$ to obtain unit Frobenius norm.
On a diagonal block, $J_2/\sqrt2$ gives a zero phase mode, while the two symmetric matrices divided by $\sqrt2$ give the diagonal positive modes.
On each external block, the commuting and anticommuting bases give two minus and two plus modes.
Substitution into the two forms yields every entry of~\eqref{eq:eig-symm} with the stated multiplicities.
Each mode combines entries with equal metric weight, so that the basis is independent of the design.
At a critical point, the metric converts the second differential into the Hessian operator.
The coefficient ratios give the aligned spectrum for every positive design, including unequal weights on repeated data blocks.

\emph{Kernel and positivity.}
Objective invariance implies that the orbit tangent lies in the Hessian kernel.
The direct calculation gives one internal zero per equal-$d$ pair for $\varrho=+1$ and two for $\varrho=-1$, with $p/2$ additional phase zeros in the skew case.
The intrinsic nullity is
$\sum_b\binom{m_b}{2}$ or
$p/2+2\sum_b\binom{m_b}{2}=\sum_bm_b^2$, respectively.
These nullities equal the orbit dimensions. Hence the Hessian kernel equals the orbit tangent at the aligned representative.
The strict ordering of $c$ makes every remaining coefficient positive.
Objective invariance transports the Hessian forms along the orbit by congruence.
For arbitrary positive weights, this preserves the kernel dimension and positivity on each metric-normal complement.
\end{namedproof}

\subsection{Proof of Corollary~\ref{cor:diagonal-block}}
\begin{namedproof}{Proof.}
The action is free and proper by the group argument in the proof of Proposition~\ref{prop:critical-manifold}.
Under the equal-block restriction, every $R\in\mathcal C_D$ commutes with $\varSigma_P^{1/2}$.
The corresponding $U$ is orthogonal and commutes with $P$.
Formula~\eqref{eq:metric-projector-form} yields $M_{XU}^\beta=M_X^\beta$ and
$g_{XU}^\beta(Z_1U,Z_2U)=g_X^\beta(Z_1,Z_2)$.
The action is therefore isometric and defines a Riemannian quotient~\citep[Thm.~9.38]{boumal2023introduction}.

Under the equal-block restriction, the map~\eqref{eq:common-reference} sends $XU$ to $WR$.
The group correspondence in the proof of Proposition~\ref{prop:critical-manifold} identifies every realization quotient with $\mathcal Q$, and the minimizing orbit becomes the single point $[W_\star]$.
At a critical point, the quotient Hessian is the horizontal restriction of the total-space Hessian~\citep[Prop.~9.45 and Ex.~9.46]{boumal2023introduction}.
The zero modes are exactly vertical, and the quotient Hessian has the positive spectrum of the total-space Hessian.
Isometry makes this spectrum constant along the minimizing orbit.
\end{namedproof}

\subsection{Proof of Theorem~\ref{thm:eigenvalue-Hf}}
\begin{namedproof}{Proof.}
Use the common quotient identified in the proof of Corollary~\ref{cor:diagonal-block} and the mode calculation in the proof of Proposition~\ref{prop:critical-manifold}.
For every admissible metric, the retained positive modes are orthogonal to the vertical space. They represent a common basis of the quotient tangent.
Restricting the two forms gives $H_0$ and $G$, with Hessian matrix $G^{-1}H_0$ and the positive ratios in~\eqref{eq:eig-symm}.
Proposition~\ref{prop:commutation} realizes every positive choice of the metric parameters. Thus the displayed metrics exhaust the family in Theorem~\ref{thm:eigenvalue-Hf}.
\end{namedproof}

\subsection{Proof of Lemma~\ref{lem:LP-invariance}}
\begin{namedproof}{Proof.}
For $\varrho=+1$, the internal log-eigenvalues $\ell^\nu_{ij}$ and external log-eigenvalues $\ell^\mu_{ij}$ both translate by $-(2\bar a+t)$.
        For $\varrho=-1$, the same calculation applies to $\ell^{\nu,\pm}_{ij}$, $\ell_i^0$, and $\ell^{\mu,\pm}_{ij}$.
        The common translation leaves the spread $\overline\lambda-\underline\lambda$ unchanged.
        Choose $t$ and $\bar a$ to impose $\tau=0$ and $x_1=0$.
        Setting $t$ to cancel $\tau$ recovers the metric-parameter redundancy of Remark~\ref{rem:rescaling}.
\end{namedproof}

\subsection{Proof of Proposition~\ref{prop:LP-wellposed}}
\begin{namedproof}{Proof.}
For any reduced design $z$, take $\overline\lambda=\max_k\ell_k(z)$ and $\underline\lambda=\min_k\ell_k(z)$.
        These bounds give a feasible point with finite objective.
        The objective is nonnegative on the feasible polyhedron. Attainment follows from the fact that a linear objective bounded below on a nonempty polyhedron attains its infimum.
        Exponentiating an optimizer gives positive design weights.
\end{namedproof}

\subsection{Proof of Theorem~\ref{thm:strong-duality}}
\begin{namedproof}{Proof.}
LP strong duality~\citep[Sec.~5.2.3]{vandenberghe2004convex} applies by Proposition~\ref{prop:LP-wellposed}.
        Each dual feasible pair gives $w=u-v\in\mathcal Z$ with $\sigma(w)\le1$.
        Conversely, take $w\in\mathcal Z$ with $\sigma(w)\le1$ and define $(w_+)_k=\max(w_k,0)$ and $(w_-)_k=\max(-w_k,0)$.
        Add the same nonnegative vector of mass $1-\sigma(w)$ to $w_+$ and $w_-$ to obtain a dual feasible pair.
        The correspondence preserves the objective value and proves the first equality. Normalizing nonzero $w$ by $\sigma(w)$ proves the second.
        Finally, $\sum_k w_k\ell_k(z)=\langle w,a\rangle$ and $\underline\lambda\le\ell_k(z)\le\overline\lambda$ give $|\langle w,a\rangle|\le\sigma(w)(\overline\lambda-\underline\lambda)$.
\end{namedproof}

\subsection{Proof of Corollary~\ref{cor:gain}}
\begin{namedproof}{Proof.}
At $(x,y,\tau)=(0,0,0)$ every log-eigenvalue equals its offset $a_k$, and $(\overline\lambda,\underline\lambda)=(\max_k a_k,\min_k a_k)$ is feasible.
        LP optimality gives $G_\star\le\max_k a_k-\min_k a_k$. Exponentiating the right-hand side gives the aligned condition number in Remark~\ref{rem:special-cases} at $(\varLambda_P,\varLambda_Q)=(I,I)$.
        For $\varrho=+1$, take $\breve Q=I_n$ and $\breve P=I_p$.
        This gives $(A,B)=(C,D)$ in~\eqref{eq:equivalent-replacement}, and~\eqref{eq:commutation-conditions} holds with $(P,Q)=(I_p,I_n)$.
        With $\beta=1$, the matrix $E=[X,X_\perp]$ is orthogonal and $O=I_{n-p}$. The metric~\eqref{eq:metric-canonical} is Euclidean.
        The identity $f=h$ and Theorem~\ref{thm:eigenvalue-Hf} give the original Hessian spectrum at the minimizer.
        For $\varrho=-1$, the condition $AQ=QA$ reads $CJ_n=J_nC$, which fails for generic SPD $C$.
\end{namedproof}

\section{Preconditioned Metrics and Their Spectra}
\label{sec:preconditioned-metric}
\label{sec:metric-costs}
We derive the spectra for the two preconditioned metrics applied to $h(Y)=\mathrm{tr}(Y^\top CYD)$ on $\mathrm{St}(p,n)$.
Write $\operatorname{sym}(A)=(A+A^\top)/2$ and
$\mathfrak L=\operatorname{sym}(Y^\top CYD)$.
For a self-adjoint ambient SPD operator $\mathcal H(Y)$, the metric is
$g_Y^{\mathcal H}(Z_1,Z_2)=\langle Z_1,\mathcal H(Y)[Z_2]\rangle$.
Its tangent projection has the form
$\mathcal P_{\mathrm T_Y}^{\mathcal H}(V)=V-\mathcal H^{-1}(YS)$, where $S=S^\top$ solves
\begin{equation}\label{eq:H-proj-S}
 \operatorname{sym}\bigl(Y^\top\mathcal H^{-1}(YS)\bigr)
 =\operatorname{sym}(Y^\top V).
\end{equation}
Positive definiteness of $\mathcal H(Y)$ makes the tangency system invertible.
The gradient and critical-point Hessian are
\begin{align}
 \mathrm{grad}\,h(Y)&=\mathcal P_{\mathrm T_Y}^{\mathcal H}\bigl(\mathcal H^{-1}(2CYD)\bigr),\\
 \mathrm{hess}\,h(Y)[Z]&=\mathcal P_{\mathrm T_Y}^{\mathcal H}\bigl(\mathcal H^{-1}(2CZD-2Z\mathfrak L)\bigr),\quad Y\in\mathcal Y.\label{eq:H-hess}
\end{align}
$\mathcal Y$ denotes the metric-independent critical set.
To verify~\eqref{eq:H-hess}, differentiate
$\mathrm{grad}\,h=\mathcal H^{-1}(2CYD-YS)$.
At a critical point the residual vanishes and $S=2\mathfrak L$.
The vanishing residual removes derivatives of $\mathcal H^{-1}$ and the connection terms.
The remaining term $\mathcal H^{-1}(Y\mathrm DS[Z])$ is normal and vanishes under projection.
At a critical point, the Rayleigh numerator is the Euclidean Stiefel Hessian form~\citep[Sec.~4.1]{absil2013extrinsic}, namely
$2\mathrm{tr}(Z^\top CZD-Z^\top Z\mathfrak L)$.
At points outside $\mathcal Y$, the diagnostic operator~\eqref{eq:H-hess} omits the metric-derivative terms of the full Hessian.

\subsection{Lagrangian Preconditioned Metric}\label{sec:metric-Lagrangian}
Following~\citet[Sec.~3.3, Case~I]{mishra2016riemannian}, set
\begin{equation}\label{eq:metric-lagrangian}
 \mathcal H_{\mathrm{Lag}}(Y)[Z]=2CZD-2\omega Z\mathfrak L,
 \qquad g_Y^{\mathrm{Lag}}(Z_1,Z_2)=\langle Z_1,\mathcal H_{\mathrm{Lag}}(Y)[Z_2]\rangle.
\end{equation}
In an eigenbasis of $C$, this operator splits into blocks
$2\lambda_a(C)D-2\omega\mathfrak L$, where $\lambda_a(C)$ ranges over the eigenvalues of $C$.
The identity $\mathrm{tr}(D^{-1}\mathfrak L)=\mathrm{tr}(Y^\top CY)>0$ implies that the largest generalized eigenvalue of $(\mathfrak L,D)$ is positive.
For $\omega\ge0$, ambient positive definiteness is equivalent to
\begin{equation}\label{eq:omega-star}
 0\le\omega<\omega_\star(Y):=\lambda_{\min}(C)/\lambda_{\max}(\mathfrak L,D).
\end{equation}
The experiments use
\begin{equation}\label{eq:omega-adaptive}
 \omega(Y)=\eta\omega_\star(Y),\qquad \eta=0.8.
\end{equation}
This pointwise rule need not be smooth when the largest generalized eigenvalue is multiple.
The inverse action solves the Sylvester equation $2CWD-2\omega W\mathfrak L=R$.

\begin{proposition}\label{prop:spectrum-Lagrangian}
For SPD data $C,D$, let $Y\in\mathrm{St}(p,n)$ be an aligned global minimizer of~\eqref{eq:quadratic-problem}.
Let $c_1\ge\cdots\ge c_n>0$ and $d_1\ge\cdots\ge d_p>0$ be the eigenvalues of $C^{-1}$ and $D$ in the minimizing order.
Under~\eqref{eq:metric-lagrangian} with $0\le\omega<\omega_\star(Y)$ from~\eqref{eq:omega-star}, the tangent Hessian spectrum is
\begin{equation*}
 \left\{\frac{(d_i-d_j)(c_i-c_j)}{c_id_i+c_jd_j-\omega(c_id_j+c_jd_i)}\right\}_{i<j\le p}
 \cup\left\{\frac{c_i-c_j}{c_i-\omega c_j}\right\}_{i\le p<j\le n}.
\end{equation*}
Each indexed value contributes one tangent mode, with multiplicities added when values coincide.
\end{proposition}
\begin{remark}\label{rmk:lag-limits}
Under Proposition~\ref{prop:spectrum-Lagrangian}, aligned frames give $\mathfrak L=\mathrm{diag}(d_i/c_i)$ and $\omega_\star=c_p/c_1$.
For simple $c$ and $p>1$, this is less than one.
The formal limit $\omega\to1^-$ sends every positive mode to one, but lies outside the ambient SPD interval.
Equal-$d$ modes remain zero.
Each positive mode is $[1+(1-\omega)r]^{-1}$ with $r>0$, making the positive condition number nonincreasing on the admissible interval.
\end{remark}

\subsection{Left--Right Preconditioned Metric}\label{sec:metric-LR}
Let $R_\varepsilon=(\mathfrak L^2+\varepsilon I_p)^{1/2}$, with $\varepsilon>0$.
The left--right approximation of~\citet[Sec.~3.2]{gao2025optimization} gives
\begin{equation}\label{eq:metric-LR}
 \mathcal H_{\mathrm{LR}}(Y)[Z]=CZR_\varepsilon,
 \qquad g_Y^{\mathrm{LR}}(Z_1,Z_2)=\mathrm{tr}(Z_1^\top CZ_2R_\varepsilon).
\end{equation}
Its inverse is $\mathcal H_{\mathrm{LR}}^{-1}[V]=C^{-1}VR_\varepsilon^{-1}$.
The tangency equation~\eqref{eq:H-proj-S} is a generalized Sylvester equation for symmetric $S$.
The experiments fix $\varepsilon=10^{-3}$.
\begin{proposition}\label{prop:spectrum-LR}
For the SPD data, aligned Stiefel minimizer, and ordered $c_i,d_i$ of Proposition~\ref{prop:spectrum-Lagrangian}, use metric~\eqref{eq:metric-LR} with $\varepsilon>0$.
Put $s_i=\sqrt{(d_i/c_i)^2+\varepsilon}$.
The tangent Hessian spectrum is
\begin{equation*}
 \left\{\frac{2(c_i-c_j)(d_i-d_j)}{c_is_i+c_js_j}\right\}_{i<j\le p}
 \cup\left\{\frac{2d_i(c_i-c_j)}{c_is_i}\right\}_{i\le p<j\le n},
\end{equation*}
with one contribution per indexed tangent mode.
\end{proposition}
\begin{namedproof}{Proof of Propositions~\ref{prop:spectrum-Lagrangian} and~\ref{prop:spectrum-LR}.}
The selected frames give $Y^\top CY=\mathrm{diag}(1/c_i)_{i\le p}$. The complementary block is
$Y_\perp^\top CY_\perp=\mathrm{diag}(1/c_j)_{j>p}$.
Internal modes are $Z_{ij}=Y(e_ie_j^\top-e_je_i^\top)$ and external modes are
$W_{ij}=Y_\perp e_{j-p}e_i^\top$.
The common Hessian-form numerators for the internal and external modes are, respectively,
$2(d_i-d_j)(c_i-c_j)/(c_ic_j)$ and $2d_i(c_i-c_j)/(c_ic_j)$.
The Lagrangian metric denominators are
$2[c_id_i+c_jd_j-\omega(c_id_j+c_jd_i)]/(c_ic_j)$ and
$2d_i(c_i-\omega c_j)/(c_ic_j)$.
The left--right denominators are
$(c_is_i+c_js_j)/(c_ic_j)$ and $s_i/c_j$.
These modes span the tangent space and simultaneously diagonalize the Hessian and metric forms.
Taking the coefficient ratios gives both spectral formulas.
\end{namedproof}

\begin{remark}\label{rmk:LR-limit}
Under Proposition~\ref{prop:spectrum-LR}, as $\varepsilon\to0$ the internal and external values tend to
$2(c_i-c_j)(d_i-d_j)/(d_i+d_j)$ and $2(c_i-c_j)$.
If $C^{-1}$ is not normalized, the external limits have no universal upper bound. They also have no uniform positive lower bound.
The data fix these weights. The constraint design provides a separate weight for each complementary canonical mode.
\end{remark}
\begin{remark}\label{rem:sp-obstruction}
The metrics~\eqref{eq:metric-lagrangian} and~\eqref{eq:metric-LR} use the symmetric multiplier of the Stiefel constraint.
A symplectic multiplier is $J_p$-self-adjoint, but its square can be nonsymmetric.
For example, the matrix with rows $(1,0,0,-3)$, $(0,1,2,0)$, $(0,3,4,0)$, $(-2,0,0,4)$ has this property for $J_p=\mathrm{blkdiag}(J_2,J_2)$.
Extending these metrics to the symplectic constraint therefore requires a separate metric construction.
\end{remark}

In our dense implementation, \textsf{PM-Lag} assembles the tangency system in $p(p+1)/2$ symmetric basis matrices and applies $n$ dense $p\times p$ solves per inverse action.
Without factorization reuse, gradient conversion costs $O(n^2p^3+np^5+p^6)$.
\textsf{PM-LR} assembles a $p^2\times p^2$ Kronecker system, giving $O(n^2p+np^2+p^6)$.
These costs apply to the direct implementations. Structured solvers and factorization reuse can reduce them.
For the metric induced by the \textsf{PC} formulation, $K=I-XP^{-1}X^\top Q$ gives
$(M_X^\beta)^{-1}=\beta^{-1}X(PP^\top)^{-1}X^\top+KK^\top$.
A structured action costs $O(n^2p+np^2+p^3)$, whereas our implementation explicitly forms $KK^\top$ and adds $O(n^3)$.
Cost evaluation, line search, and retraction are additional operations.

\section{Technical Lemmas}
\label{appendix:lemmas}

    The notation for scalars and dimensions in the following lemmas is local to this section. All eigenvalue multiplicities are algebraic.

    \begin{lemma}\label{lem:eigenvalue-skew}
        Let $p,q$ be positive integers and $a_i,c_i$ $(1\le i\le p)$, $b_j,d_j$ $(1\le j\le q)$ be real scalars.
        Let $T\in\mathbb{R}^{4p^2\times 4p^2}$ be the permutation matrix satisfying $\mathrm{vec}(\varOmega^\top)=T\mathrm{vec}(\varOmega)$ for each $\varOmega\in\mathbb{R}^{2p\times 2p}$.
        Let
        \begin{equation*}
            \begin{aligned}
                \varXi &= \mathrm{diag}(a_1,\ldots,a_p) \otimes I_2 \otimes \mathrm{diag}(b_1,\ldots,b_q) \otimes I_2 \\
                &\quad - \mathrm{diag}(c_1,\ldots,c_p) \otimes J_2 \otimes \mathrm{diag}(d_1,\ldots,d_q) \otimes J_2.
            \end{aligned}
        \end{equation*}
        Then the following statements hold.
        \begin{enumerate}
            \item The eigenvalues of $\varXi$ are $a_i b_j \pm c_i d_j$, each with multiplicity $2$, for $1\le i\le p$ and $1\le j\le q$.
            \item If $p=q$, the eigenvalues of $(I+T)\varXi$ are $p(2p-1)$ zeros together with
                \begin{enumerate}[label=(\roman*)]
                    \item $(a_i b_j + a_j b_i) \pm (c_i d_j + c_j d_i)$, each with multiplicity $2$, for $1\le i<j\le p$.
                    \item $2(a_i b_i + c_i d_i)$, with multiplicity $2$, for $1\le i\le p$.
                    \item $2(a_i b_i - c_i d_i)$, with multiplicity $1$, for $1\le i\le p$.
                \end{enumerate}
                The three families contain $p(2p+1)$ values with the stated multiplicities.
                Each vanishing listed value adds to the zero eigenvalue multiplicity.
        \end{enumerate}
    \end{lemma}

    \begin{namedproof}{Proof.}
        Choose independent eigenvectors $v_1,v_2$ of $J_2$ with eigenvalues $\eta_1=\mathrm i$ and $\eta_2=-\mathrm i$, where $\mathrm i^2=-1$.
        The tensor vectors $w_{i,s,j,t}=e_i\otimes v_s\otimes e_j\otimes v_t$ form an eigenbasis of $\varXi$, with eigenvalues $\mu_{i,s,j,t}=a_ib_j-c_id_j\eta_s\eta_t$. The choices of $s,t$ give (a).
        For $p=q$, $T$ exchanges $w_{i,s,j,t}$ and $w_{j,t,i,s}$.
        A distinct transpose pair with eigenvalues $\mu,\nu$ gives the matrix $\left[\begin{smallmatrix}\mu&\nu\\\mu&\nu\end{smallmatrix}\right]$, with eigenvalues $0$ and $\mu+\nu$.
        Each fixed coordinate gives $2\mu$.
        For $i<j$, the two choices with $s=t$ give the plus family and the two with $s\ne t$ give the minus family.
        For $i=j$, the two fixed coordinates give $2(a_ib_i+c_id_i)$ and the remaining pair gives $2(a_ib_i-c_id_i)$.
        The $\binom{2p}{2}=p(2p-1)$ distinct transpose pairs account for all redundant zeros.
    \end{namedproof}

    \begin{lemma}\label{lem:eigenvalue-symm}
        Let $p,q$ be positive integers and $a_i,c_i$ $(1\le i\le p)$, $b_j,d_j$ $(1\le j\le q)$ be real scalars.
        Let $T\in\mathbb{R}^{p^2\times p^2}$ denote the permutation matrix such that $\mathrm{vec}(\varOmega^\top)=T\mathrm{vec}(\varOmega)$ for each $\varOmega\in\mathbb{R}^{p\times p}$.
        Let
        \begin{equation*}
            \varXi = \mathrm{diag}(a_1,\ldots,a_p) \otimes \mathrm{diag}(b_1,\ldots,b_q) - \mathrm{diag}(c_1,\ldots,c_p) \otimes \mathrm{diag}(d_1,\ldots,d_q).
        \end{equation*}
        Then the following statements hold.
        \begin{enumerate}
            \item The eigenvalues of $\varXi$ are $a_i b_j - c_i d_j$, each with multiplicity $1$, for $1\leq i\leq p$ and $1\leq j\leq q$.
            \item If $p=q$, the eigenvalues of $(I-T)\varXi$ are $p(p+1)/2$ zeros together with all values of $(a_i b_j + a_jb_i) - (c_i d_j + c_j d_i)$, each with multiplicity $1$, for $1\leq i<j\leq p$.
        \end{enumerate}
    \end{lemma}

    \begin{namedproof}{Proof.}
        Each $e_i\otimes e_j$ has eigenvalue $\mu_{ij}=a_ib_j-c_id_j$ under $\varXi$.
        When $p=q$, $T$ exchanges $e_i\otimes e_j$ and $e_j\otimes e_i$.
        On each pair with $i<j$, $(I-T)\varXi$ has matrix $\left[\begin{smallmatrix}\mu_{ij}&-\mu_{ji}\\-\mu_{ij}&\mu_{ji}\end{smallmatrix}\right]$, with eigenvalues $0$ and $\mu_{ij}+\mu_{ji}$.
        Together with the $p$ vanishing diagonal coordinates, these pairs give $p(p+1)/2$ zeros.
    \end{namedproof}
\end{document}